\documentclass[twoside,12pt]{article}

\usepackage{amsmath}
\usepackage{amssymb}
\usepackage{pxfonts}
\usepackage{graphicx}
\usepackage{eucal}
\usepackage{mathrsfs}
\usepackage{theorem}
\usepackage{pifont}
\usepackage{amscd}
\usepackage{enumitem}
\usepackage{float}

\usepackage{subfig}
\usepackage{caption}

\usepackage{sectsty}
\usepackage{pseudocode}
\usepackage{color}
\usepackage{fancyhdr}
\usepackage{framed}
\usepackage[all]{xy}
\usepackage[pagebackref]{hyperref}

\definecolor{shadecolor}{rgb}{0.8,0.8,0.8}

\theoremheaderfont{\fontfamily{pzc}\bfseries\large}

\newtheorem{theorem}{Theorem}[section]
\newtheorem{lemma}[theorem]{Lemma}
\newtheorem{proposition}[theorem]{Proposition}
\newtheorem{corollary}[theorem]{Corollary}
\newtheorem{definition}[theorem]{Definition}

{\theorembodyfont{\rmfamily}

\newenvironment{proof}{{\flushleft \emph{Proof}:}}{\hfill\ding{110}}
\newenvironment{comment}{{\flushleft \fontfamily{pzc}\bfseries\large Comment:}}{}
\newenvironment{comments}{{\flushleft \fontfamily{pzc}\bfseries\large Comments:}}{}

\newcommand{\dashint}{\dashint}

\newcommand{\secref}[1]{Section~\ref{#1}}

\newcommand{\thmref}[1]{Theorem~\ref{#1}}
\newcommand{\defref}[1]{Definition~\ref{#1}}

\newcommand{\propref}[1]{Proposition~\ref{#1}}
\newcommand{\lemref}[1]{Lemma~\ref{#1}}
\newcommand{\corref}[1]{Corollary~\ref{#1}}

\newcommand{\tr}{\operatorname{tr}}
\newcommand{\Cof}{\operatorname{Cof}}

\renewcommand{\ker}{\operatorname{Ker}}
\newcommand{\dist}{\operatorname{dist}}

\newcommand{\trg}{\operatorname{tr}_{\g}}
\newcommand{\weakly}[1]{\stackrel{#1}{\rightharpoonup}}

\newcommand{\br}{V_{\N}/V_{\M}}

\newcommand{\bdx}{\partial_x}
\newcommand{\bdy}{\partial_y}

\newcommand{\supp}{\operatorname{supp}}
\newcommand{\Lipp}{\operatorname{Lip}_+}

\newcommand{\Lippi}{\operatorname{Lip}_{+}^{\operatorname{inj}}}

\newcommand{\R}{\mathbb{R}}
\newcommand{\D}{\mathcal{D}}
\newcommand{\Sig}{\Sigma}
\newcommand{\sig}{\sigma}

\newcommand{\diff}{\operatorname{Diff}}

\newcommand{\SO}{\operatorname{SO}}

\newcommand{\Otwom}{\operatorname{O}^{-}_2}

\newcommand{\id}{\operatorname{Id}}
\renewcommand{\div}{\operatorname{div}}
\newcommand{\Vol}{\operatorname{Vol}}
\newcommand{\Volg}{\operatorname{Vol}_\g}

\newcommand{\g}{\mathfrak{g}}
\newcommand{\h}{\mathfrak{h}}

\newcommand{\covder}{\nabla}

\newcommand{\M}{\mathcal{M}}

\newcommand{\N}{\mathcal{N}}

\newcommand{\til}{\tilde}

\newcommand{\brk}[1]{\left(#1\right)}          
\newcommand{\sAverage}[1]{\langle#1\rangle}      

\newcommand{\deriv}[2]{\frac{d#1}{d#2}}
\newcommand{\pd}[2]{\frac{\partial#1}{\partial#2}}

\newcommand{\beq}{\begin{equation}}
\newcommand{\eeq}{\end{equation}}

\newcommand{\qrt}{\frac{1}{4}}
\newcommand{\bqrt}{1/4}

\providecommand{\half}{\frac{1}{2}}

\newcommand{\distSO}[1]{\dist(#1,\SO)}

\newcommand{\IP}[2]{\sAverage{#1,#2}}

\newcommand{\al}{\alpha}
\newcommand{\be}{\beta}

\newcommand{\de}{\delta}
\newcommand{\ep}{\epsilon}

\newcommand{\diag}{\operatorname{diag}}
\newcommand{\pl}{\partial}

\newcommand{\SOd}{\operatorname{SO}_d}
\newcommand{\SOn}{\operatorname{SO}_n}

\newcommand{\SOtwo}{\operatorname{SO}_2}
\newcommand{\CO}{\operatorname{CO}_2}

\newcommand{\GLp}{\operatorname{GL}_n^+}
\newcommand{\GLtwo}{\operatorname{GL}_2^+}

\newcommand{\lam}{\lambda}

\renewcommand{\skew}{\operatorname{skew}}

\def\Xint#1{\mathchoice
{\XXint\displaystyle\textstyle{#1}}%
{\XXint\textstyle\scriptstyle{#1}}%
{\XXint\scriptstyle\scriptscriptstyle{#1}}%
{\XXint\scriptscriptstyle\scriptscriptstyle{#1}}%
\!\int}
\def\XXint#1#2#3{{\setbox0=\hbox{$#1{#2#3}{\int}$ }
\vcenter{\hbox{$#2#3$ }}\kern-.6\wd0}}

\renewcommand{\dashint}{\Xint-}

\newcommand{\ab}[1]{{\color{black} #1}}

\numberwithin{equation}{section}

\begin{document}
\title{Embedding surfaces inside small domains with minimal distortion}
\author{Asaf Shachar\footnote{Institute of Mathematics, The Hebrew University of Jerusalem.} }
\date{}
\maketitle

\begin{abstract}
Given two-dimensional Riemannian manifolds $\M,\N$, we prove a lower bound on the distortion of embeddings $\M \to \N$, in terms of the areas' discrepancy $V_{\N}/V_{\M}$, for a certain class of distortion functionals. For $V_{\N}/V_{\M} \ge 1/4$, homotheties, provided they exist, are the unique energy minimizing maps attaining the bound, while for $V_{\N}/V_{\M} \le 1/4$, there are non-homothetic minimizers. We characterize the maps attaining the bound, and construct explicit non-homothetic minimizers between disks. We then prove stability results for the two regimes. We end by analyzing other families of distortion functionals. In particular we characterize a family of functionals where no phase transition in the minimizers occurs; homotheties are the energy minimizers for all values of $V_{\N}/V_{\M}$, provided they exist. 
\end{abstract}
\tableofcontents
\section{Introduction}
\subsection{Setting}
Let $\M, \N$ be connected, compact, oriented smooth $2$-dimensional Riemannian manifolds (possibly with Lipschitz boundaries) having areas $V_{\M}, V_{\N}$. Suppose that $V_{\N} \le V_{\M}$. We consider the following general question:
\begin{quote}
How to embed $\M$ in $\N$ with minimal distortion? 
\end{quote}
Let $\Lippi(\M,\N)$ be the space of injective almost everywhere Lipschitz maps $\M \to \N$ having almost everywhere non-negative Jacobian. 
By injective a.e.~we mean that $|f^{-1}(q)| \le1$ for a.e.~$q \in \N$.
 We consider various \emph{distortion} functionals $E:\Lippi(\M,\N) \to [0,\infty)$ 
 and provide lower bounds on $E$ in terms of the discrepancy between $V_{\M},V_{\N}$.  Intuitively, such a bound must exist since squeezing a domain into a smaller domain must carry distortion. 

%

The functional $E$ is required to satisfy $E(\phi)=0$ if and only if $\phi$ is an orientation-preserving isometric immersion.
 We assume that   
$E(\phi)=\int_{\M} W(d\phi)\Volg$,
where $W \ge 0$ is some non-negative energy density and $\Volg$ is the Riemannian volume form of $\M$. Set 
\beq
\label{eq:gen_min}
E_{\M,\N}=\inf_{\phi \in \Lippi(\M,\N)} E(\phi).
\eeq

This minimization problem is motivated by the theory of incompatible elasticity, which is a branch of elasticity concerned with bodies that do not have a stress-free reference configuration (see, e.g., \cite{ESK08,KES07,KVS11, DCGLL13,AKMMS16}). Such bodies are typically modelled as Riemannian manifolds $(\M, \g)$, with the ambient space being another manifold $(\N, \h)$.  
The physical context of the present problem is that of ``confinement" where the elastic body is constrained within some ambient environment. 
%

As a first example, consider the case where $\M=\Omega \subseteq \R^2, \N=\R^2$. Then $\phi:\Omega \to \R^2$ is a map between flat spaces and $W:M_2 \to \R$, where $M_2$ is the space of real $2 \times 2$ matrices. We assume that $W$ is bi-$\SOtwo$ invariant, i.e., $W(RA)=W(AR)=W(A)$ for every $R \in \SOtwo$; this implies that $W$ is a function of the singular values of its argument. 
%
For the most part of this work, we assume the prototypical energy density $W_p(A):=\dist^p(A,\SOtwo)$ where $\dist(.,\SOtwo)$ is the Euclidean
distance from $\SOtwo$. This energy density can be defined similarly for mappings between Riemannian manifolds:
Denote by $\g,\h$ the metrics on $\M,\N$ respectively.
$\SOtwo$ is replaced by $\SO(\g_x,\h_y)$---the set of orientation preserving isometric linear maps $T_x\M\to T_y\N$. Given 
$\phi:\M \to \N,x\in\M$, set $\SO(\g,\phi^*\h)_x:=\SO(\g_x,\h_{\phi(x)})$; the Riemannian metrics on $\M,\N$ induce an inner-product on $\operatorname{Hom}(T_x\M,T_{\phi(x)}\N)$. We measure the distance of $d\phi_x$ from $\SO(\g,\phi^*\h)_x$ using the distance induced by this inner product. 
Define
\beq
\label{eq:p_energy}
E_p(\phi)=\dashint_{\M} W_p(d\phi)\Volg, \, \, \, \text{ where } \, \,  W_p(d\phi)=\dist^p\brk{d\phi,\SO(\g,\phi^*\h)}, 
\eeq
$\dashint$ denotes the integral divided by the volume of $\M$. For convenience,  we may write $\SOtwo$ instead of $\SO(\g,\phi^*\h)$, even when referring to the manifold case.

A natural question is whether $E_{\M,\N}$ is attained, and if it does, to characterize the energy minimizing maps. 
Since $W_p(A)=\dist^p(A,\SOtwo)$ is not quasiconvex (see \cite{Sil01} or \cite{Dol12}), it is not clear a-priori that minimizers exist. 

Note that non-injective maps may have lower energy than injective maps: e.g.~an isometric immersion from a circle of radius $2$ into a circle of radius $1$ has zero energy, even though there is a discrepancy between the lengths. Thus length discrepancy alone cannot be used to construct a lower bound on the energy of non-injective maps. Higher dimensional examples can be obtained analogously.

A particular case of interest is where $\N$ is a smaller scaled copy of $\M$. Recall that $\phi:(\M,\g) \to (\N,\h)$ 
is called a homothety if $\phi^*\h=\lam^2 \g$ for some constant $\lam >0$. Equivalently:
\[
\h_{\phi(p)}(d\phi_p(v),d\phi_p(w))=\lam^2 \g_{p}\brk{v,w} \, \, \,  \text{ for every }  \, \, p \in \M   \, \, \,  \text{ and } \, \, v, w \in T_p\M. 
\]
We say that $\M,\N$ are \emph{homothetic} if there exists an orientation-preserving diffeomorphism $\phi:\M \to \N$ which is a homothety. If $\phi^*\h=\ab{\lam^2} \g$, we denote $\M=\lam \N$. \ab{(The length scale was multiplied by $\lam$).} Note that every Lipschitz homothety is smooth: If $\phi \in W^{1,\infty}(\M;\N)$ satisfies $d\phi\in \lam \SO(\g,\phi^*\h)$ a.e.~then 
 $\phi \in C^{\infty}$ (see e.g.~ \cite{Har58,CH70,Tay06,kupferman2019reshetnyak}). 
Two homothetic surfaces $\M,\N=\lam \M$  have the same geometry up to the scale $\lam$. 
One might conjecture that in that case the minimizers of $E_p$ are the homotheties. Surprisingly, this is not always the case. 

Notation: From this point forward, we omit $\Volg$ from all integrals, i.e.~we denote $\int_{\M} f \Volg$  by $\int_{\M} f$. 
%


\subsection{Energy bounds and exact minimizers}
\label{sec:Energy bounds and exact minimizers}


To formulate our main theorem, we introduce the following notation: Given $A \in M_2$, we denote by $\sig_1(A),\sig_2(A)$ its singular values. Set 
\[
K= \{ A \in M_2 \, | \, \det A \ge 0 ,  \, \, \sig_1(A)+\sig_2(A)=1\}.
\] 
\label{def:doublewell_combined}
By the AM-GM inequality $\sig_1(A)+\sig_2(A) =1$ implies $\det A \le 1/4$;  for $0 \le s \le 1/4$ set 
\[
K_s=\{ A \in K \, | \, \det A=s\}. 
\]
Singular values are defined for linear maps between inner product spaces; for $\phi:\M \to \N$, we write $d\phi \in K$ a.e.~if $d\phi_p \in K$ for almost every $p \in \M$; $d\phi_p$ is a map $T_p\M \to T_{\phi(p)}\N$, so the metrics on $\M,\N$ are implicitly involved in the definition. 


Define $F:[0, \infty) \to \mathbb R$ by
\beq
\label{def:energy_point_vol_bd}
F(s) =
\begin{cases}
1-2s, & \text{ if  }\, 0 \le s \le \frac{1}{4} \\
2(\sqrt{s}-1)^2,  & \text{ if  }\, s \ge \frac{1}{4} 
\end{cases}
\eeq


Our main results are the following: 
\begin{theorem}
\label{the:main_Euclid_bound_injective}
Let $p \ge 2$. For every $\phi \in \Lippi(\M,\N)$
\beq
\label{eq:volboundgen8a}
E_p(\phi) \ge  F^{p/2}\brk{\frac{V_{\phi\brk{\M} } }{V_{\M}}  }.
\eeq
If $V_{\phi\brk{\M} }/{V_{\M}}  \ge \bqrt $ equality holds  if and only if $\phi$ is a homothety. 

If $V_{\phi\brk{\M} }/{V_{\M}}  \le \bqrt$ and $p=2$, equality holds if and only if $d\phi \in K$ a.e.

If $V_{\phi\brk{\M} }/{V_{\M}}  \le \bqrt$ and $p>2$, equality holds if and only if $d\phi \in K_{V_{\phi\brk{\M}} / V_{\M} }$ a.e.

%
%
\end{theorem}

We have required $p \ge 2$, since our analysis relies on the convexity of $F^{p/2}$, which is not valid for $p<2$. We expect the convex envelope of $F^{p/2}$ to play a role in the analysis for $1 \le p <2$, which we leave for future works.
\begin{corollary}
\label{cor:surject_min}
If $V_{\N} \le V_{\M}$, then $\,   E_p(\phi) \ge  F^{p/2}\brk{ \br  }$, with equality if and only if $\phi$ is surjective and satisfies the conditions above. In particular, if there exists a bijection with the required properties, it is energy-minimizing.
\end{corollary}
When squeezing a body into a smaller environment, it might seem profitable to use all the space given. This heuristic does not always apply, however---one can add a very thin ``neck" to $\N$ while barely changing its area; there is no reason for an optimal embedding to fill in this neck. It is therefore an interesting question to characterize when there exists a surjective minimizer, and when all minimizers are surjective. 

\begin{corollary}
\label{cor:homothet_min}
If $\M,\N$ are homothetic, and $\br \ge \bqrt$, then the homotheties are the unique energy minimizers. 
\end{corollary}

Let $\M,\N$ be homothetic,  $\br \ge \bqrt$; by \corref{cor:homothet_min} there exists a unique energy minimizing diffeomorphism up to a composition with an isometry. We shall see that this is not always the case when $\br < \bqrt$. 

Suppose that $\br < \bqrt$. For $p>2$ the minimizers attaining the bound $F^{p/2}\brk{\br}$ lie in $K_s$. 
 The well $K_s$ is flexible---it contains many non-affine smooth maps, see \cite{351550}. 
Thus the minimizers in the regime $\br < \bqrt$ have more flexibility compared to the rigid homothetic case when $\br \ge \bqrt$.

A natural question is whether the bound $F^{p/2}(V_{\N}/V_{\M})$ is attained for every $\M,\N$ satisfying $V_{\N}/V_{\M} < 1/4$. 
Among smooth maps this is not always the case---there is a topological obstruction. If we take $\N=\lam \M$, $\M$ a closed surface, then the existence of a map $\phi \in C^1(\M,\N)$ with $d\phi \in K_{\br}$ implies that $\M$ is diffeomorphic to a torus, see \cite{RobObs}. Moreover, for some metrics on the torus, a \emph{discretization} of the admissible singular values (which in turn corresponds to a discretization of the admissible compressions ratios $\lambda$) may happen; see e.g.~\cite{375931} for the flat torus.
\paragraph{Concrete example}
Consider the case when $\M, \N=\lam \M$ are disks and $0<\lam \le 1/2$. We prove the following: 
\begin{proposition}
\label{prop:min1disks}
Let $0<\lam \le 1/2$, and let $\D \subseteq \mathbb{R}^2$ be the closed unit disk. Denote by $ \D^o=\D\setminus\{0\}$ the disk with the origin removed. Then for any $p\ge 2$,
\[
E_{\D^o , \lam \D^o}=\min_{\phi \in \Lippi(\D^o , \lam \D^o)} E_p(\phi)=F^{p/2}\brk{\lam^2  }
\] 
is realized by a smooth diffeomorphism.
\end{proposition}

Let $\diff(\M,\N)$ be the set of orientation-preserving smooth diffeomorphisms $\M \to \N$. Through approximation one deduces from \propref{prop:min1disks} that
\beq
\label{eq:inf_disks}
E_{\D,\lam \D}=\inf_{\phi \in \diff(\D , \lam \D)} E_p(\phi)=F^{p/2}\brk{\lam^2  }.
\eeq
Note that this is a statement about complete disks.
It is an interesting question whether this infimum is attained. The following proposition answers it affirmatively for $p=2$. 
\begin{proposition}
\label{prop:min2disks}
Let $0<\lam <1/2$, and let $\D \subseteq \mathbb{R}^2$ be as above. Then
\[
\min_{\phi \in \diff(\D , \lam \D)} E_2(\phi)=F\brk{\lam^2  },
\]
and there exists an infinite-dimensional family of energy-minimizing diffeomorphisms. 
\end{proposition}
Let $\lam < 1/2$; for $p>2$ minimizers have constant singular values, whereas for $p=2$ only the sum of their singular values is constant.  This additional freedom enables us to construct smooth minimizers between complete disks for $p=2$; we do not know whether this is possible for $p>2$. We also do not know whether the minimizers $\D^o \to \D^o$ attaining the bound $F^{p/2}\brk{ \frac{V_{\N}}{V_{\M}}  }$ for $p>2$ are unique. (In the proof of \propref{prop:min1disks} we construct one minimizer for each value of $p$.)





\paragraph{Symmetry breaking and connection to physics}
For $\lam \ge 1/2$ the minimizers $\D\to \lam \D$ are radially-symmetric, given by $\brk{r,\theta } \mapsto \brk{\psi(r),\theta }$, whereas  for $\lam < 1/2$ there are no radially-symmetric minimizers (see \secref{sec:no_radial_min}); at a certain threshold of compression, the radially-symmetric maps stop being minimizers. This symmetry breaking resembles physical phenomena observed in metamaterials under compression. 

One example---a ``holes experiment"---is described in \cite[figure 5]{bertoldi2017flexible}. Another example was demonstrated experimentally in \cite{stoop2015curvature}; As a result of applying isotropic pressure on a ball (a polydimethylsiloxane-coated elastomer), a wrinkling pattern on the boundary occurs. 
These phenomena are a byproduct of the competition between two energy terms---stretching and bending. 



%
%

In contrast, the analysis in the current work suggests that a "bulk" symmetry breaking may occur as a byproduct of pressure with stretching energy alone (no bending). It is an interesting question whether such a phase transition can be observed experimentally. The transition presented here occurs at compression ratio $\lam=1/2$, which might seem unrealistic, since after such a large compression the material would no longer remain elastic. The precise value of $1/2$, however, is model-dependent; it might be possible that for some materials, the ratio is closer to $1$ and thus more realistic (see also the comment below \thmref{thm:bijective_homo_min_short}.)

\paragraph{The well $K$ and critical points}

\thmref{the:main_Euclid_bound_injective}  singles out maps $\phi \in \Lippi(\M,\N)$ with $d\phi \in K$ as energy minimizing maps.
It is therefore natural to wonder whether there is any direct connection between the well $K$ and critical points of the energy $E_p$.
The next result clarifies this:
\begin{proposition}
\label{prop:wellK_is_critic}
Let $\phi \in C^2(\M,\N)$ with $J\phi>0$ and $d\phi \in K$. Then $\phi$ is a critical point of $E_2$, and it is a critical point of $E_p$ for $p \neq 2$ if and only if its singular values are constant.
\end{proposition}
$d\phi \in K$ implies $J\phi \ge 0$; we required here $J\phi>0$ since the integrand of $E_p$ is smooth only when restricted to invertible matrices.

Another context where the well $K$ arises is the following:
\begin{proposition}
\label{prop:constant_sing}
Let $\Omega \subseteq \R^2$ be an open connected domain, and let $\phi \in C^2(\Omega,\R^2)$ have constant singular values. Suppose that $\phi$ is $E_p$ critical for some $p \ge 1$. Then $\phi$ is either affine or satisfies $d\phi \in K$. 
\end{proposition}
Due to \propref{prop:constant_sing}, when looking for critical maps between Euclidean spaces having constant singular values, one is naturally led to study maps in the well $K$. We do not know whether  \propref{prop:constant_sing} holds for maps between arbitrary surfaces.

We prove propositions \ref{prop:wellK_is_critic} and \ref{prop:constant_sing} by establishing an alternative characterization of the well $K$ (see \propref{prop:geometric_cof_char_well}). We deduce from this characterization that the solutions to the Euler-Lagrange equation of $E_p$ do not have to be $C^1$, even in the Euclidean case.

\subsection{Rigidity}

A natural question is the nature of minimizing sequences, $\phi_n 
\in \Lippi(\M,\N)$ satisfying $E_p(\phi_n) \to \inf E_p$. We treat the case of $p=2$, and distinguish between the cases of $V_{\N}/V_{\M} < 
1/4$ and $V_{\N}/V_{\M} >1/4$; the flexibility of minimizing sequences is quite different between these two cases.

\subsubsection{Rigidity for $\bqrt < \br \le 1$ }

Assume that  $\bqrt < \br \le 1$. 
By \thmref{the:main_Euclid_bound_injective}, $ E_2(\phi) \ge  F( \br)$ and equality holds if and only if $\phi$ is a bijective homothety. The following result is an asymptotic version of \thmref{the:main_Euclid_bound_injective}:

\begin{theorem}
\label{thm:asymptot_hom}
Let $\M,\N \subseteq \mathbb{R}^2$ be open, bounded sets with Lipschitz boundaries.
Suppose that  $1/4 < V_{\N}/V_{\M} \le 1$. Let $\phi_n \in \Lippi(\M,\N)$ and assume that $ E_2(\phi_n) \to  F(\br)$. Then $\phi_n$ has a subsequence converging strongly in $W^{1,2}(\M,\N)$ to a smooth surjective homothety $\phi \in \Lippi(\M,\N)$ with Jacobian $J\phi=\br$, which is injective on $\M^\circ$ and satisfies $\phi(\M^\circ) \subseteq \N^\circ$. If $\phi_n(\pl \M)\subset \pl \N$, then $\phi$ is injective on $\M$, and if $\partial \M, \partial \N$ are smooth, then $\phi$ is smooth up to the boundary, and is a bijective diffeomorphic homothety.
\end{theorem}

The assumption $\phi_n(\pl \M)\subset \pl \N$ cannot be dropped.  
Take for example $\M=[-1,1]^2$, and let $\N=\M/\sim$ be the flat $2$-torus with $\sim$ the standard equivalence relation.
Then $\phi_n:\M\to \N$ given by $\phi_n(x)=(1-1/n)x$ are injective and satisfy all the other conditions, but converge uniformly to the quotient map $\pi:\M\to \N$, which is obviously not an isometry but merely an isometric immersion.

\thmref{thm:asymptot_hom} states that if the infimal energy is that of a homothety, then any minimizing sequence converges to a homothety. In particular, \emph{there exists} a homothety between $\M $ and $\N$. (Note that we do not assume a-priori that $\M$, $\N$ are homothetic.)
This result is analogous to a classical result of Reshetnyak \cite{Res67}: 
\begin{quote}
\emph{
Let $\Omega\subset \R^d$ be an open, connected, bounded domain, $1\le p<\infty$.
If $\phi_n\in W^{1,p}(\Omega;\R^d)$ satisfy $\int_\Omega \phi_n\,dx = 0$ and $\dist(d \phi_n, \SOd)\to 0$ in $L^p(\Omega)$, then $\phi_n$ has a subsequence converging strongly in $W^{1,p}(\Omega;\R^d)$ to an isometric mapping.
}
\end{quote}
Reshetnyak's theorem was generalized to mappings between manifolds in \cite{kupferman2019reshetnyak}.
Reshetnyak's theorem states that a sequence of mappings whose $p$-energy tends to that of an isometry, converges (modulo a subsequence) to an isometric immersion. \thmref{thm:asymptot_hom} is the analogous result obtained by replacing "isometry" with "homothety".

Note that we restricted \thmref{thm:asymptot_hom} to Euclidean domains. Most of the proof holds as is for arbitrary surfaces; however, there is a key element which holds only for Euclidean domains; it is the so-called ``higher integrability property of determinants", which states that if
$\phi_n \weakly{} \phi$ in $W^{1,2}(\M,\N)$ and $J\phi_n \ge 0$, then $ J\phi_n \rightharpoonup J\phi $ in $L^1(K)$ for any compact $K  \Subset \M^\circ$, see \cite{Mul90}.
This statement does not hold between manifolds; for example, there is a sequence of conformal diffeomorphisms of the sphere $\mathbb{S}^2$, which converges in $W^{1,2}$ to a constant. (see e.g.~\cite[p. 415]{giaquinta1998cartesian1} , or \cite{381194}). Generalizing  \thmref{thm:asymptot_hom} to general surfaces is an interesting problem left for future works. (For reasons of brevity we treated here stability only for the case $p=2$; we expect a similar result should hold for $ p>2$.)

In \thmref{thm:asymptot_hom} we required  $ V_{\N}/V_{\M} >1/4$.
When $ V_{\N}/V_{\M} =1/4$ there may be minimizing sequences which do not converge to homotheties: e.g., when $\M,\N$ are Euclidean disks, take any minimizer for the problem $\M \to \lam_n \M$, where $\lam_n \nearrow 1/2$ and scale it by $1/(2\lam_n)$. 
\subsubsection{Rigidity for $\br \le \bqrt$}

Assume that $V_{\phi\brk{\M} }/{V_{\M}} \le \bqrt$. By \thmref{the:main_Euclid_bound_injective}, $ E_2(\phi) \ge  F( V_{\phi\brk{\M} }/{V_{\M}} )$ and equality holds if and only if $d\phi \in K$. Similarly, if $\br \le \bqrt$, then by \corref{cor:surject_min}, $E_2(\phi) \ge  F( \br )$ and equality holds if and only if $d\phi \in K$ and $\phi$ is surjective. The following rigidity estimate is a quantitative sharpening of these statements. 


\begin{theorem}
\label{thm:low_regime_rigidity}
Let $\phi \in \Lippi(\M,\N)$. Suppose that $V_{\phi\brk{\M} }/{V_{\M}}  \le \bqrt$. Then 
\beq
\label{eq:quantbound_vol_distortion_dist_SO2_1abc}
\dashint_{\M} \dist^2(d\phi,K)  \le E_2(\phi) - F\brk{\frac{V_{\phi\brk{\M} }}{V_{\M}}}  \le  2\dashint_{\M} \dist^2(d\phi,K) .
\eeq
Furthermore, if $\br \le \bqrt$ then 
\beq
\label{eq:quantbound_vol_distortion_dist_SO2_1abcde}
\dashint_{\M}  \dist^2(d\phi,K)+\frac{2(V_{\N}-V_{\phi(\M)}) }{V_{\M}}   \le E_2(\phi) - F\brk{\frac{V_{\N}}{V_{\M}}} \le 2 \dashint_{\M} \dist^2(d\phi,K)  +\frac{2(V_{\N}-V_{\phi(\M)}) }{V_{\M}}.
\eeq
Thus $\lim_{n \to \infty}E_2(\phi_n) = F(\br)$ if and only if $V_{ \phi_n\brk{\M} } \to V_{\N}$ and $\dist(d\phi_n,K)$ converges to $0$ in $L^2$.
\end{theorem}
Equation \eqref{eq:quantbound_vol_distortion_dist_SO2_1abc} has an interesting interpretation: 
Minimizing $E_2$ roughly means ``get as close as you can to the $\SOtwo$-well"; Equation  \eqref{eq:quantbound_vol_distortion_dist_SO2_1abc} implies that when adding a constraint on the areas ($V_{\phi\brk{\M} }/V_{\M} \le \bqrt$), and trying to approach the value $F(V_{\phi\brk{\M} }/V_{\M} )$, we get an equivalent problem of getting close to a different well---the well $K$. So, 
\[
(\SOtwo\text{-well problem}) + \text{(area constraint)} \simeq (K\text{-well problem}).
\]


\subsection{Other distortion functionals}

\label{sec:Other_distortion_functionals}
A natural question is whether the occurrence of phase transitions depends on the functional $E$. 
We demonstrate that indeed, it depends crucially on properties of the energy density $W$. Here is the setting: 

Let  $f:(0,\infty) \to [0,\infty)$ be a continuous function satisfying $f(1)=0$, which is strictly increasing on $[1,\infty)$, and strictly decreasing on $(0,1]$. We think of $f$ as a cost function measuring how much $x$ deviates from $1$. Every such $f$ induces a functional $E_f:\Lipp(\M,\N) \to [0,\infty]$ by
\beq
\label{eq:gen_energy_functional}
E_f(\phi)=\dashint_{\M}  f(\sigma_1(d\phi))+ f(\sigma_2(d\phi)) .
\eeq

The energy density $W=f(\sigma_1)+f(\sigma_2)$ incorporates as a special case Ogden-like materials whose energy density is given by
$f(x)=\Sigma_{k=1}^N a_k (x^{\alpha_k}-1)$ (see e.g.~\cite[p. 189]{Cia88}). The energy $E_2$ is recovered as a special case by setting $f(x)=(x-1)^2$. 
This setting does not cover the case of $E_p$ for $p>2$; however, it covers the classical $p=2$ case as well as other natural examples (see \secref{subsec:logarithmic_example} ). It is not hard to generalise the analysis to functionals of the form 
\[
E_g(\phi)=\int_{\M} g(\sigma_1(d\phi),\sigma_2(d\phi)),
\]
where $g$ does not necessarily decompose into an additive sum of contributions from the singular values (this is the most general form of a bi-$\SOtwo$ invariant density.) 

For the class of energies \eqref{eq:gen_energy_functional}, we have two main results. The first is that homotheties are always energy minimizing when compressing by a small amount: If $\lam$ is sufficiently close to $1$, then
the energy minimizing maps $\M \to \lam \M$ are the homotheties. 

Our second result singles out a large family of cost functions $f$, for which no phase transition occurs; the homotheties, if exist, remain the energy minimizers for any degree of compression.


To state our results, we define the auxiliary function
\beq
\label{eq:min_uniform0}
F_f(s)=\min_{xy=s,x,y>0} f(x)+ f(y), \, \,  \, \, \text{for } \, \, s \in (0,\infty).
\eeq
When  $f(x)=(x-1)^2$ this definition of $F_f$ agrees with Definition \eqref{def:energy_point_vol_bd}.

A function which is $k$-times differentiable at $x_0$ is called \emph{flat} at $x_0$ if all its derivatives vanish at $x_0$. If a function is not flat at a local minimum, then it is strictly convex in some neighbourhood of that minimum. 

Our first result is the following:

\begin{theorem}
\label{thm:bijective_homo_min_short} 
Suppose that $f$ is differentiable and not flat at $x=1$. Then $\exists \al \in (0,1)$ such that for every $\M,\N$ satisfying $\al \le \br \le 1$, and every $\phi \in \Lippi(\M,\N)$,
\[
 E_f(\phi) \ge F_f\brk{\frac{V_{\N}}{V_{\M}}},
\]
with equality if and only if $\phi$ is a surjective homothety. In particular, if $\M,\N$ are homothetic the homotheties are the unique energy minimizers. 

Moreover, if $V_{\phi(\M)}/{V_{\M}} \ge \al$, then $E_f(\phi) \ge  F_f\brk{V_{\phi(\M)}/{V_{\M}} },$ with equality if and only if $\phi$ is a homothety.
\end{theorem}
\begin{comment}
$\al$ depends on $f$; replacing $f$ with $\til f(x)=f(x^r)$ results in $\til F(s)=F(s^r)$, so the transition point where the homotheties stop being minimizers can be pushed arbitrarily close to $1$.
\end{comment}


Our second result is the following:

\begin{theorem}
\label{thm:inject_gen_bound}
Let $g:\R \to [0,\infty)$ be a continuous function, which is strictly decreasing and strictly convex on $(-\infty,0]$, and strictly increasing on $[0,\infty)$, with $g(0)=0$. Set $f(x) = g(\log x)$ and let $E_f,F_f$ be as in equations \eqref{eq:gen_energy_functional} and \eqref{eq:min_uniform0}. Then for any $\phi \in \Lippi(\M,\N)$
\beq
\label{eq:volboundgen4bca}
E_f(\phi) \ge  F_f\brk{\frac{V_{\phi(\M)}}{V_{\M}}  },
\eeq
and equality holds  if and only if $\phi$ is a homothety. Thus, if there exists an injective homothety, it is energy-minimizing among the maps whose images have the same volume.
Finally, if $V_{\N} \le V_{\M}$ then $\,   E(\phi) \ge  F_f\brk{\br}$, with equality if and only if $\phi$ is a surjective homothety. 
\end{theorem}

\paragraph{Related works}
The works \cite{sivaloganathan2010global} and \cite{mora2010explicit} study minimal distortion under uniaxial extension for incompressible materials ($J\phi=1$ is assumed); in \cite{mora2010explicit} the authors show that a phase transition in the energy minimizers occurs when the expansion becomes sufficiently large.



\paragraph{Structure of this paper}

In \secref{sec:Euclidean_distortion_functional} we prove the lower bound \thmref{the:main_Euclid_bound_injective} and Propositions \ref{prop:min1disks}, \ref{prop:min2disks} on exact minimizers between disks. In \secref{sec:approx_high1} we prove the stability Theorems \ref{thm:asymptot_hom} and \ref{thm:low_regime_rigidity}. In \secref{sec:secgeom_prop_well} we prove \propref{prop:wellK_is_critic} and \propref{prop:constant_sing}. In \secref{sec:gen_distortion_functional_prf} we prove Theorems \ref{thm:bijective_homo_min_short} and \ref{thm:inject_gen_bound} regarding general distortion functionals. In \secref{sec:Discuss} we discuss some open questions that arise from this work.

\section{Volume bounds for the Euclidean functional}
\label{sec:Euclidean_distortion_functional}
\subsection{Pointwise bound}
We begin with the following lower bound on  $\dist(A,\SOtwo)$ in terms of $\det A$: 
\begin{lemma}
 \label{lem:bound_vol_distortion_dist_SO22}
 
Let $A \in M_2$ satisfy $\det A \ge 0$. Then 
\[
 \dist^2(A,\SOtwo)   \ge F(\det A)=\begin{cases}
 1-2\det A, & \text{ if  }\, 0 \le \det A \le \frac{1}{4} \\
2(\sqrt{\det A}-1)^2,  & \text{ if  }\, \det A \ge \frac{1}{4}. 
\end{cases}
\]

If $\det A \ge  \bqrt$ equality holds  if and only if $A$ is conformal, i.e.~$\sig_1(A)=\sig_2(A)$.

If $\det A \le \bqrt$ equality holds if and only if $\sig_1(A)+\sig_2(A)=1$, or $A \in K$, where $K$ is defined in \ref{sec:Energy bounds and exact minimizers}.



\end{lemma}

 \lemref{lem:bound_vol_distortion_dist_SO22} can be proved separately for the cases where $\det A \ge 1/4$ and $\det A \le 1/4$. We will state and prove quantitative generalisations of it for both regimes. The claim for $\det A \ge 1/4$ follows from \lemref{lem:quantbound_vol_distortion_dist_confSO22}, and the claim for $\det A \le 1/4$ follows from \lemref{lem:quantbound_vol_distortion_dist_SO22} (see comment after Equation \eqref{eq:quantbound_vol_distortion_dist_SO2_2}.)

We give a direct proof of  \lemref{lem:bound_vol_distortion_dist_SO22} (which does not rely on the results just stated) in Appendix \ref{sec:add_proof}.
%

Naively one might expect that conformal matrices are the closest to $\SOtwo$ in the class of matrices with a given determinant. Here is a heuristic argument why this is false when the determinant is sufficiently small:  Putting equal sharing of the distortion on the singular values is suboptimal in the marginal case of $s=0$; setting the singular values to be $(0,1)$ is better than the conformal option which is $(0,0)$. The same heuristic works when $s\ll1$-it is profitable to set one singular value close to zero and the other one close to $1$.
\subsection{Proof of \thmref{the:main_Euclid_bound_injective}}

We prove a more general result, which does not assume injectivity. We denote by 
 $\Lipp(\M,\N)$ the space of Lipschitz maps $\phi:\M \to \N$ with $J\phi \ge 0$ a.e.~

\begin{theorem}
\label{the:main_Euclid_bound}
Let $p \ge 2$ and let $\phi \in \Lipp(\M,\N)$. Then
\beq
\label{eq:volboundge}
E_p(\phi) \ge  F^{p/2}\brk{\dashint_{\M}  J\phi  }. 
\eeq
If $\dashint_{\M} J\phi \ge 1/4 $ equality holds if and only if $\phi$ is a homothety. 

If $\dashint_{\M} J\phi \le 1/4 $ and $p=2$, equality holds if and only if $d\phi \in K$ a.e.~

If $\dashint_{\M} J\phi \le 1/4 $ and $p>2$, equality holds if $J\phi$ is constant a.e.~and $d\phi \in K_{J\phi}$ a.e.
\end{theorem}
 \thmref{the:main_Euclid_bound_injective} follows immediately from \thmref{the:main_Euclid_bound}:
For $\phi \in \Lippi(\M,\N)$ the area formula \cite{hajlasz1999sobolev} implies that $\int_{\M}J\phi =\int_{\N} |f^{-1}(q)| = V_{\phi\brk{\M} } $. 

For the proof we shall need the following auxiliary result:  
\begin{lemma}
\label{lem:Jensen_equality1}
Let $X$ be a probability space, and let $g:X \to \mathbb [0,\infty) $ be in $L^1(X)$. 
Let $F:\mathbb [0,\infty) \to [0,\infty)$ be convex and strictly convex on $[a,\infty)$ for some  $a \in (0,\infty)$.

Then, $\int_X g  \in [a,\infty)$ implies that $\int_X F\circ g=F(\int_X g)$ if and only if $g$ is constant a.e.

$\int_X g  \in (0,a]$ implies that if $\int_X F\circ g=F(\int_X g)$, then $g \le a$ a.e.
\end{lemma}

A proof is given in appendix \ref{sec:convexity_results}.


\begin{proof}[of \thmref{the:main_Euclid_bound}] We begin with the case $p=2$. The function $F$ in Definition \eqref{def:energy_point_vol_bd} is convex: A direct computation shows that $F \in C^1$, and 
\beq
\label{def:energy_point_vol_bd_deriv}
F'(s)=\begin{cases}
 -2, & \text{ if  }\, 0 \le s\le\frac{1}{4} \\
 2\left(1-\frac1{\sqrt{s}}\right), & \text{ if  }\, s\geq\frac14,
\end{cases} 
\eeq
is continuous and non-decreasing. By \lemref{lem:bound_vol_distortion_dist_SO22}
\beq
\label{eq:Euclideanbound1}
E_2(\phi)= \dashint_{\M}  W_2(d\phi)  \stackrel{(1)}{\ge} \dashint_{\M}  F(J\phi)  \stackrel{(2)}{\ge} F\brk{\dashint_{\M}  J\phi  },
\eeq
where inequality $(2)$ follows from the convexity of $F$.


  Suppose that $\dashint_{\M}J\phi \ge \bqrt$:
Using the strict convexity of $F$ on $[\frac14,\infty ]$, and applying \lemref{lem:Jensen_equality1} with $g=J\phi,a=\bqrt$, we deduce that inequality $(2)$ is an equality if and only if $J\phi$ is constant a.e., which implies $J\phi \ge 1/4$ a.e.~By  \lemref{lem:bound_vol_distortion_dist_SO22} inequality $(1)$ is an equality if and only if $d\phi$ is conformal a.e.~Thus $E_2(\phi) =  F(\dashint_{\M}  J\phi )$ if and only if $\phi$ is a homothety.

Suppose that $\dashint_{\M}J\phi \le \bqrt$:
Again by  \lemref{lem:Jensen_equality1}, inequality $(2)$ is an equality if and only if $J\phi \le \bqrt$ a.e.~\lemref{lem:bound_vol_distortion_dist_SO22} implies that inequality $(1)$ is an equality if and only if $d\phi \in K$ a.e.~This completes the proof for $p=2$.

Suppose that $p>2$, and set $q=p/2$. $F^q$ is strictly convex. 
Since $x \to x^q$ is monotonic, $W_2(A) \ge F(\det A) \Rightarrow W_p(A)=W_2^q(A) \ge F^q(\det A)$. Modifying Equation \eqref{eq:Euclideanbound1} we get
%
\beq
\label{eq:Euclideanbound_p}
E_p(\phi)= \dashint_{\M}  W_p(d\phi)   \stackrel{(1)}{\ge} \dashint_{\M}  F^q(J\phi)  \stackrel{(2)}{\ge} F^q\brk{\dashint_{\M}  J\phi  }
\eeq
with inequality $(2)$ being an equality if and only if $J\phi$ is constant a.e.~

The remaining step of the proof---characterizing the equality case in inequality $(1)$ is exactly as for $p=2$.  
\end{proof}
\begin{proof}[of \corref{cor:surject_min}]
Since
$
0 \le \frac{V_{\phi\brk{\M} }}{V_{\M}} \le \frac{V_{\N}}{V_{\M}} \le1
$
and $F|_{[0,1]}$ is strictly decreasing, \thmref{the:main_Euclid_bound_injective} implies that
\beq
\label{eq:Euclideanbound_pbc}
E_p(\phi) \ge  F^{p/2}\brk{\frac{V_{\phi\brk{\M} } }{V_{\M}}  }  \stackrel{(3)}{\ge} F^{p/2}\brk{\frac{V_{\N}}{V_{\M}}},
\eeq
with inequality $(3)$ being an equality if and only if $V_{\phi\brk{\M} } = V_{\N}$. Since $\phi$ is continuous, $\phi(\M)$ is compact, in particular closed in $\N$. $V_{\phi\brk{\M} } = V_{\N}$ implies that $\phi(\M)$ is dense, hence it coincides with $\N$.
\end{proof}
\subsection{Non-homothetic minimizers}
\label{Non_homothetic_min}
In this section we prove Propositions \ref{prop:min1disks} and \ref{prop:min2disks} regarding the existence of energy minimizing diffeomorphisms between disks.


We say that $ \phi_1,  \phi_2 :\M \to \N$ are \emph{equivalent} if there exist isometries $T:\M \to \M, S:\N \to \N$ such that
$ \phi_1=S  \circ \phi_2 \circ T$.
Since $W$ is bi-$\SOtwo$ invariant, two equivalent maps have the same energy, so the set of energy minimizers is invariant under compositions of isometries, and thus forms a union of equivalence classes. 

We introduce the following notation; for $0 \le \sig_1 \le \sig_2$, set 
\[
K_{\sig_1,\sig_2}=\{ A \in M_2 \, | \, \det A \ge 0 ,  \, \, \sig_1(A)=\sig_1,\sig_2(A)=\sig_2\}.
\]

\subsubsection{Reducing the homothetic problem to self-maps}
We want to prove statements about energy minimizers between homothetic disks. We formulate a more general result which characterizes energy minimizers $\M \to \lam \M$ for arbitrary surfaces $\M$, not necessarily disks.

Let $0<\lam \le 1/2$. We prove a correspondence between diffeomorphic minimizers $\M \to \lam \M$ attaining the bound $F^{p/2}(\lam^2)$ and certain diffeomorphisms $\M \to \M$ whose singular values have specific properties. 




Fix a homothety $\psi \in \diff(\M,\lam\M)$. 
Let $0<\sigma_1<\sigma_2$ be the unique numbers satisfying $\sigma_1+\sigma_2=1,   \sigma_1 \sigma_2=\lambda^2$, and let $\phi \in \Lippi(\M,\lam\M)$. 

Since 
\[
d\phi \in K_{\lam^2} \iff d\phi \in K_{\sig_1,\sig_2} \iff d(\psi^{-1} \circ \phi) \in K_{ \frac{\sig_1}{\lam}, \frac{\sig_2}{\lam} },
\] the map $T: \phi \mapsto \psi^{-1} \circ \phi$ is a bijection 
 \[
 \{ \phi:\M \to \lam\M \,\, |\, \,d\phi \in K_{\lam^2}\} \to   \{ \til \phi:\M \to \M \,\, |\, \, d\til \phi \in K_{\frac{\sig_1}{\lam}, \frac{\sig_2}{\lam} } \}.
 \]
Similarly, $T$ maps bijectively 
\[
\{ \phi:\M \to \lam\M \, \, |  \, \, d\phi \in K \} \, \, \text{ to } \, \,\{ \til \phi:\M \to \M \, | \, \,  \sig_1(d\til \phi)+\sig_2(d\til \phi)=\frac{1}{\lam} \}.
\] 
Since $ \phi_1,  \phi_2 :\M \to \lam\M$  are equivalent  if and only if $T(\phi_1), T(\phi_2)$ are equivalent, $T$ induces a map between the corresponding equivalence classes.

Let $p>2$; by \thmref{the:main_Euclid_bound_injective}, the value $F^{p/2}(\lam^2)$ is attained by $\phi \in \diff(\M , \lam\M)$  if and only if  $d\phi \in  K_{\lam^2}$. For $p=2$, the energy value $F(\lam^2)$ is attained by $\phi \in \diff(\M , \lam\M)$  if and only if  $d\phi \in  K$. This implies the following: 
\begin{theorem}
\label{thm:equiv_min_self_map}
Let $0<\lam \le 1/2$. Then for $p>2$. 
\[
\min_{\phi \in \diff(\M , \lam\M)} E_p(\phi)=F^{p/2}(\lam^2)
\]
if and only if there exists an \textbf{area-preserving} diffeomorphism $\phi:\M \to \M$ having constant sum of singular values $1/\lam$.
The equivalence classes of energy minimizing diffeomorphisms $\M \to \lam\M$ are then isomorphic (through $T$) to equivalence classes of \textbf{area-preserving} diffeomorphisms $\M \to \M$  whose singular values sum up to $1/\lam$. 

Similarly, for $p=2$ we have
\[
\min_{\phi \in \diff(\M , \lam\M)} E_2(\phi)=F(\lam^2) 
\]
if and only if there exists a diffeomorphism $\phi:\M \to \M$ having constant sum of singular values $1/\lam$. The equivalence classes of energy minimizing diffeomorphisms $\M \to \lam\M$ are isomorphic to equivalence classes of diffeomorphisms $\M \to \M$  whose singular values sum up to $1/\lam$. 

%
\end{theorem}
A similar statement can be formulated for bijective Lipschitz minimizers instead of diffeomorphic minimizers.
\begin{comments}
\begin{itemize}
\item The property of being area-preserving with constant singular values is inverse-invariant: If
 $\phi \in \diff(\M)$ is area-preserving with singular values $\sig_1,\sig_2$ then so is $\phi^{-1}$ (this is special to dimension $2$).
Since the energy minimizers are in one-to-one correspondence with such diffeomorphisms, this suggests an approach for finding non-equivalent minimizers. However, $\phi, \phi^{-1}$ may already be equivalent, so we do not always obtain new minimizers in this way.
\item
Let $\lam=1/2$ and $\phi \in \diff(\M)$.  Then,
\[
\label{com:com_crit_case_rigid}
\sig_1(d\phi)+\sig_2(d\phi)=\frac{1}{\lam}=2 \Rightarrow \sqrt{J\phi} \le \half \big(\sig_1(d\phi)+\sig_2(d\phi)\big) =1.
\]
Thus, $1=\dashint_{\M} J\phi \le \dashint_{\M} 1=1$, which implies $J\phi=1$; thus we have an AM-GM \emph{equality}, which implies $\sig_1(d\phi)=\sig_2(d\phi)=1$. So for $\lam=1/2$, all the
self-diffeomorphisms mentioned in \thmref{thm:equiv_min_self_map} are isometries. This recovers the rigid case, where the energy minimizers are homotheties.
\end{itemize}
\end{comments}

\subsubsection{Minimizers having constant singular values}

We prove \propref{prop:min1disks} regarding the existence of energy minimizing diffeomorphisms between punctured disks.
 By \thmref{thm:equiv_min_self_map} it suffices to construct an area-preserving diffeomorphism $\phi \in \diff( \D^o)$ having constant singular values $\sigma_1\le \sigma_2$ satisfying $\sigma_1+\sigma_2= 1/ \lam $.  When $\lam=1/2$, the only such diffeomorphisms are isometries (see comment after \thmref{thm:equiv_min_self_map}.)
 

Let $0<\lam < 1/2$. 
The linear map $x \to \begin{pmatrix} \sigma_1 & 0 \\\ 0 & \sigma_2 \end{pmatrix}x \notin \diff( \D^o)$ since its image gets outside of $\D^o$ as $ \sigma_2 > 1$, thus we must consider non-affine candidates. 
We construct $\phi \in \diff( \D^o)$ with $d\phi \in K_{\sig_1,\sig_2}$:

Let $c \in \R$, and define $\phi_c: \D^o \to  \D^o$ in polar coordinates by
\beq
\label{eq:logarithmic_min}
\phi_c:\big(r,\theta\big )\mapsto \big(r,\theta+c\log(r)\big).
\eeq
$\phi_c$ is the flow over time $c$ of the divergence free vector field $\log r \frac{\partial}{\partial \theta}$.

The figures below describe the action of $\phi_c$ for $c=0,1/2$. For $c=0$, $\phi_0=\id$; figure $(a)$ describes the disk, divided into coloured pizza-slices. $\phi_{1/2}$ maps each pizza-slice into a twisted slice, as can be seen in figure $(b)$. The slices (twisted or not) all have equal area.



\begin{figure}[H]
    \centering
    \subfloat[$c=0$]{{\includegraphics[width=3.5cm]{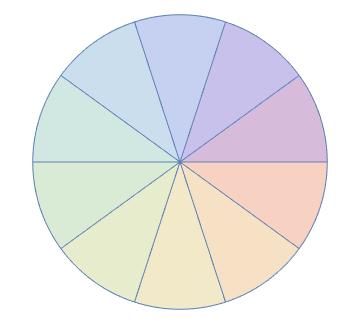} }}
    \qquad
    \subfloat[$c=1/2$]{{\includegraphics[width=3.5cm]{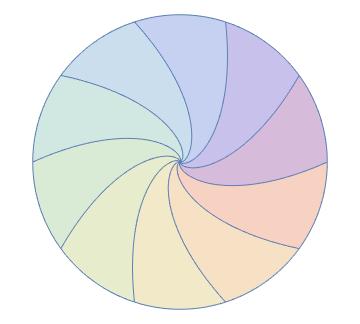} }}
\end{figure}



For the sake of more general examples, it is convenient to analyze more general maps of the form 
\beq
\label{eq:concentric_min}
\phi:\big(r,\theta\big )\mapsto \big(\psi(r),\theta+h(r)\big).
\eeq
We have
\[ d\phi\brk{\frac{\partial}{\partial \theta}(r,\theta)}=\frac{\partial}{\partial \theta}\big(\phi(r,\theta)\big), d\phi\brk{\frac{\partial}{\partial r}(r,\theta)}=\psi'(r)\frac{\partial}{\partial r}\brk{\phi(r,\theta)}+h'(r)\frac{\partial}{\partial \theta}\brk{\phi(r,\theta)},\]
so w.r.t the orthonormal frame $\{ \frac{\partial}{\partial r},\frac{1}{r}\frac{\partial}{\partial \theta}\}$, 
\[
d\phi\brk{\frac{\partial}{\partial r}(r,\theta)}=\psi'(r)\frac{\partial}{\partial r}\brk{\phi(r,\theta)}+\big(h'(r)\psi(r) \big) \bigg(\frac{1}{\psi(r)} \cdot \frac{\partial}{\partial \theta}\brk{\phi(r,\theta)}\bigg),
\]
and 
\[ 
d\phi\brk{\frac{1}{r}\frac{\partial}{\partial \theta}(r,\theta)}=\frac{\psi(r)}{r}\cdot \brk{\frac{1}{\psi(r)}\frac{\partial}{\partial \theta}\brk{\phi(r,\theta)}}.
\]
In other words, 
\beq
\label{eq:disks_example1}
[d\phi]_{\{ \frac{\partial}{\partial r},\frac{1}{r}\frac{\partial}{\partial \theta}\}}=\begin{pmatrix} \psi' & 0 \\\ h'\psi & \frac{\psi}{r}\end{pmatrix}.
\eeq
Specializing to $\phi_c$, where $\psi(r)=r, h(r)=c\log(r)$, we get
\beq
[d\phi_c]_{\{ \frac{\partial}{\partial r},\frac{1}{r}\frac{\partial}{\partial \theta}\}}=\begin{pmatrix} 1 & 0 \\\ c & 1\end{pmatrix}.
\eeq

Set $A_c=\begin{pmatrix} 1 & 0 \\\ c & 1\end{pmatrix}$, and let $\sig_1(c) \le \sig_2(c)$ be its singular values. 
Since $\sig_1(c)  \sig_2(c)=1$ it follows that $0 < \sig_1(c)  \le 1$.
Since 
\[
|A_c|^2=2+c^2=\sig_1^2(c) +\sig_2^2(c) \le 2\sig_2^2(c),
\]
\[
\lim_{c \to \infty}\sig_2(c)=\infty \, \, \text{  and } \, \, \lim_{c \to \infty}\sig_1(c)=0. 
\]
On the other hand $\lim_{c \to 0}\sig_1(c)=\sig_1(0)=1$. By continuity, $\sig_1(c)$ attains all the values in $(0,1]$ when $c$ ranges over $\R$.

Next, we observe that $\sig_i(c)=\sig_i(-c)$: this follows from $\phi_c^{-1}=\phi_{-c}$ which holds since $\phi_c$ is a flow (this is also immediate from the definition of $\phi_c$). In fact $\sig_i(c)=\sig_i(\tilde c)$ if and only if $c=\pm \tilde c$. Indeed, $\sig_i(c)$ are uniquely determined by the sum $\sig_1^2(c) +\sig_2^2(c)=2+c^2$ and the product $\sig_1^2(c)  \sig_2^2(c)=1$. Thus $\sig_i(c)=\sig_i(\tilde c)$ if and only if $c^2=\til c^2$. 

Thus, for every $\sig_1 \in (0,1)$, there exists a unique $c>0$ such that $\sig_1=\sig_1(\pm c)$; every possible value of $\sig_1$ is attained exactly twice when $c$ ranges over $\R$, except for $\sig_1=1$ which is attained only once for $c=0$. To conclude, we proved the following:
\begin{theorem}
\label{thm:disks_exmp1}
Let $0<\sigma_1<\sigma_2$ satisfy $\sigma_1\sigma_2=1$. There exists a diffeomorphism $\phi:\D^o \to \D^o$ with $d\phi \in K_{\sig_1,\sig_2}$. 
\end{theorem}

\propref{prop:min1disks} now follows as an immediate corollary from combining \thmref{thm:disks_exmp1}
together with \thmref{thm:equiv_min_self_map}.

The mapping $\phi_c$ is not differentiable at the origin. However, it can be approximated in $W^{1,2}$ by diffeomorphisms $\D \to \D$, which proves assertion \eqref{eq:inf_disks}. The approximation can be done by interpolating between a logarithm and a constant in the phase, so the approximating map is a rotation near the origin.

Note that $\phi_c$ and $\phi_c^{-1}$ are equivalent, as $\phi_{c}^{-1}=\phi_{-c}=J\phi_{c}J$, where $J=\begin{pmatrix} 1& 0 \\ 0 & -1 \end{pmatrix}$ is the reflection around the $x$ axis. Thus we do not get non-equivalent minimizers (see comment after \thmref{thm:equiv_min_self_map}). 

\subsubsection{Minimizers having constant sum of singular values}

We prove \propref{prop:min2disks} regarding the existence of energy minimizing diffeomorphisms between complete disks.


Let $0<\lam < 1/2$, and set $\al=1/\lam$. By \thmref{thm:equiv_min_self_map} we need to construct $\phi\in \diff(\D)$ having constant sum of singular values $\al$. 

For a matrix $A=\begin{pmatrix} a & 0 \\\ b & c\end{pmatrix}$ with positive determinant, 
\[
\sigma_1(A)+\sigma_2(A)=\al \, \, \text{ if and only if } \, \, \al^2=|A|^2+2\det A,
\]
so
\[
\sigma_1(A)+\sigma_2(A)=\al \, \, \text{ if and only if } \, \, (a+c)^2+b^2=\al^2.
\]
Equation \eqref{eq:disks_example1} implies that for maps $\phi:\D \to \D$ given by \eqref{eq:concentric_min}, the equality $ \sigma_1(d\phi)+\sigma_2(d\phi)=\al $ reduces to the ODE
\beq
\label{eq:disks_example2}
\brk{\psi'+\frac{\psi}{r}}^2+(h'\psi)^2=\al^2. 
\eeq
The following proposition asserts that for a very large family of functions $\psi$, we can find suitable functions $h$ solving Equation \eqref{eq:disks_example2}, and that $(\psi,h)$ give rise to a diffeomorphism. 

\begin{proposition}
\label{prop:many_minimizers_peq2}
Let $\al >2$. Let $\psi:[0,1] \to [0,1]$ be any concave, smooth strictly increasing function satisfying $\psi(0) = 0$,  $\psi(1) = 1$ and $ 2\psi'(0)=\al $, which is linear near zero. 
There exists a unique (up to signs and additive constants) smooth function $h:[0,1] \to \mathbb R$ such that $h$ and $\psi$ solve equation \eqref{eq:disks_example2}. The associated map $\phi$ given by  \eqref{eq:concentric_min} is then a smooth diffeomorphism of $\D$.
\end{proposition}

Explicitly, the uniqueness in $h$ is as follows: If $h_1,h_2$ solve equation \eqref{eq:disks_example2} with the same $\psi$, then either $h_1=h_2+c$ or $h_1=-h_2+c$, for some constant $c$.
This freedom in $h$ is expected: Denoting by $\phi_h$ the map given in \eqref{eq:concentric_min}, $J\phi_hJ=\phi_{-h}$, where $J$ is the reflection around the $x$ axis. Thus $\phi_h,\phi_{-h}$ are equivalent (adding a constant in $h$ amounts to composing with a rotation).
\begin{proof}
Define $g:[0,1] \to \mathbb R$ by
\beq
g(r) =
\begin{cases}
\alpha, & \text{ if  }\, r=0 \\
\psi'(r)+\frac{\psi(r)}{r},  & \text{ if  }\, 0 < r \le 1.
\end{cases}
\eeq
$g$ is non-increasing: In the neighbourhood of zero where $\psi$ is linear, $g=\al$. 
For $r>0$, 
\[
g'(r)=\psi''(r)+\frac{1}{r}(\psi'(r)-\frac{\psi(r)}{r}),
\]
and both summands are non-positive; the concavity of $\psi$ implies that $\psi'' \le 0$ and 
\[
\psi(r)=\int_0^r \psi'(t)dt \ge \int_0^r \psi'(r)dt=r\psi'(r).
\]
We prove that $g$ becomes smaller than $\al$ at some point; since it is non-increasing it remains smaller from that point onward. 

Set $t_0=\sup\{ \,r \, | \, g(r)=\al \}$.  Since $g$ is non-increasing and $g \le \al $, $g|_{[0,t_0]}=\al$. Indeed, $g^{-1}(\{\al\}) \subseteq [0,1]$ being compact implies $g(t_0)=\al$, and for $r \in [0,t_0]$ we have $\al \ge g(r) \ge g(t_0) $.


 $g|_{[0,t_0]}=\al$ implies that $\psi|_{[0,t_0]}$ is linear: $g'(r)=0 \Rightarrow \psi''(r)=0$. The assumptions  
 \[
 \psi'(0)=\frac{\al}{2}>1,\psi(0)=0, \psi(1)=1
 \] 
imply that $\psi$ cannot be linear all the way up to $r=1$,  thus $t_0<1$.
 
To conclude, we showed that there exists $t_0 \in (0,1)$ such that $g|_{[0,t_0]}=\al$ and $g|_{(t_0,1]}<\al$.  Equation \eqref{eq:disks_example2} implies that 
\[
\psi(r)h'(r)=\pm \sqrt{\al^2-g^2(r)}.
\]
We focus on the positive branch, and choose
\beq
\label{eq:disks_example3}
\psi(r)h'(r)=\sqrt{f(r)},
\eeq
where we define $f(r):=\al^2-g^2(r)$, so $f|_{[0,t_0]}=0$ and $f|_{(t_0,1]}>0$.  
This implies that $h'|_{[0,t_0)}=0$ is smooth. Similarly, $h'=\sqrt{f}/\psi$ is smooth on $(t_0,1]$, since $f$ and $\psi$ are both positive there. The only possible point of non-smoothness that we need to worry about is $t_0$. The smoothness of $h'$ at $t_0$ is a non-trivial result, that follows from the following lemma:


%
\begin{lemma}
\label{lem:smoot_root}
Let $f:[0,\epsilon) \to [0,\infty)$ be a smooth function which is strictly positive on $(0,\epsilon]$ and satisfying $f(0)=0$, $f^{(k)}(0)=0$ for every natural $k$. Then $\sqrt f$ is infinitely (right) differentiable at $x=0$, with all its right derivatives zero.
\end{lemma}

\begin{comment}
The only issue here is the infinite differentiability of $\sqrt f$; once this is established, it immediately follows that all the derivatives vanish: $f^{(k)}(0)=0$ implies $f=o(x^k)$ for every natural $k$, so $\sqrt f=o(x^k)$ as well. 
\end{comment}
The only proof of \lemref{lem:smoot_root} that we are aware of is deducing it from a sharper result. It follows as a special case from Theorem 2.2. in \cite[p. 639]{bony2010square}. (The relevant definition that is used in this theorem is definition 1.1 on p. 636.)

Applying this lemma to $f=\al^2-g^2$ (at $t_0$ instead of at $0$),  we deduce that $\sqrt f $ is infinitely differentiable at $t_0$, and hence  smooth on $[0,1]$. Since $\psi(t_0)>0$, it follows that $h'=\sqrt{f}/\psi$ is smooth at $t_0$. This completes the proof that $h'$ is smooth on $[0,1]$.  
So, given $\psi$, we have a smooth $h:[0,1] \to \mathbb R$ such that $\psi,h$ solve \eqref{eq:disks_example2}. 
We verify that the corresponding map $\phi$ given by  \eqref{eq:concentric_min} is a smooth diffeomorphism $\D \to \D$. $\phi$ is clearly bijective (since $\psi:[0,1] \to [0,1]$ is bijective), and smooth everywhere except possibly at the origin. By our construction $h|_{[0,t_0]}=c$ is constant, and $\psi(r)=\frac{\al}{2}r$ on $[0,t_0]$. Thus the restriction of $\phi$ to $\{ z \, | \, |z| \le t_0 \}$  is given by 
\[
\phi:\big(r,\theta\big )\mapsto \big(\frac{\al}{2}r,\theta+c\big),
\]
which is a simple dilation composed with a fixed rotation. (This is the reason for taking $\psi$ linear near zero, to evade the possibility of exploding phase $h$ at the origin.) 
We proved that $\phi:\D \to \D$ is a smooth bijective map. Its inverse map, given by
\[
\phi^{-1}:\big(r,\theta\big )\mapsto \big(\psi^{-1}(r),\theta-h(\psi^{-1}(r))\big)
\]
is also smooth, so $\phi$ is a diffeomorphism, as required.
\end{proof}

%
Since for a fixed $\al >2$, we can choose $\psi:[0,1] \to [0,1]$ that satisfies the conditions of \propref{prop:many_minimizers_peq2} rather arbitrarily, we constructed an infinite-dimensional family of minimizers. We note that all the minimizers we constructed between disks are of the form \ref{eq:concentric_min}. We do not know whether there exist minimizers which are not of this form.

\subsubsection{Non-existence of radial minimizers }
\label{sec:no_radial_min}
As mentioned at the end of \secref{sec:Energy bounds and exact minimizers}, when $\lam < 1/2$ there are no radially-symmetric minimizers $\D \to \lam \D$.
Here we mention two approaches for showing this:

The first is to prove that the uniform contraction $x \to \lam x$ is energy-minimizing among the radial maps. This can be done by direct computation. Thus there exists a radial minimizer if and only if homotheties are minimizing, and we already showed that this is not the case when $\lam < 1/2$. 

Alternatively, we can use our characterisation of minimizers: Any minimizer must have constant sum of singular values. Considering equations  \eqref{eq:concentric_min} and \eqref{eq:disks_example1}, we see that if $\phi(r,\theta)=(\psi(r),\theta)$ then $\sig_i(d\phi)$ are $\psi', \frac{\psi}{r}$. If their sum is constant, then $\psi$ is linear and we again get a homothety.
\label{sec:homogen_scal_best}
\section{Quantitative analysis}
\label{sec:approx_high1}

\subsection{Approximate minimizers when $1 \ge \frac{V_{\N}}{V_{\M}} > \bqrt$}
In this section we prove \thmref{thm:asymptot_hom}. As most of the steps in the proof hold for arbitrary surfaces, we will prove them in this greater generality, and assume $\M,\N$ are Euclidean domains only when it is necessary.

Let $\phi_n \in \Lippi(\M,\N)$ and assume that $ E_2(\phi_n) \to  F\brk{ \frac{V_{\N}}{V_{\M}}}$; by passing to a subsequence we may assume that $ \phi_n \rightharpoonup \phi $ in $W^{1,2}$. We need to prove that $\phi$ is a homothety; our strategy is to prove first that it is conformal, then show it has a constant Jacobian.

Since the proof is long, we divide it into several steps. First, we show that $\dashint_{\M}  J\phi_n \to \frac{V_{\N}}{V_{\M}}$ (\lemref{lem:asymptotic_Jacob_avg}); in particular $J\phi_n > \frac{1}{4}$ \emph{on average}, i.e.~$\dashint_{\M}  J\phi_n > \qrt$. Next, we prove that $J\phi_n$ is \emph{asymptotically} greater than $ \qrt$, i.e.~$\lim_{n \to \infty} V_{\{J\phi_n \le \frac{1}{4}\}}=0$ (\propref{prop:concentrated_Jacobian}). Then we show that $   \int_{\{J\phi_n \le \frac{1}{4}\}}  |d\phi_n|^2 \to 0$, i.e.~the norm of $d\phi_n$ is concentrated at $\{J\phi_n> \qrt\}$ (\lemref{lem:concentration_ophi_norm}).

The next step (\propref{prop:asymptotic_conformal}) is proving that $\phi_n$ is asymptotically conformal, i.e.~ 
\[
 \int_{\M}  \dist^2(d\phi_n,\CO) \to 0,
 \] 
and that the weak limit of an asymptotically conformal sequence is conformal (\lemref{lem:limi_asymptot_conf_is_conf}). This shows that $\phi$ is conformal. In \lemref{lem:limi_asymptot_conf_is_conf} we use the higher integrability property of determinants, thus we assume that $\M,\N \subseteq \mathbb{R}^2$. 
 The remaining step is to prove that $J\phi$ is constant. 
Together with \lemref{lem:limi_asymptot_conf_is_conf}, these are the only places where we use the assumption $\M,\N \subseteq \mathbb{R}^2$.
 
  Finally, we prove that $\phi_n$ strongly converges to $\phi$, and that the limit of an asymptotically surjective maps is surjective and injective a.e.~ (\lemref{lem:injeclimitinject}).

To motivate the different steps in the proof we begin by decomposing Inequality \eqref{eq:Euclideanbound1} into more refined steps. Denote
$A_n=\{ p \in \M \, | \, J\phi_n(p) \le \frac{1}{4} \}$, $B_n=\M \setminus A_n$; then
%
%
%
\beq
\begin{split}
\label{eq:Euclideanquant1}
E_2(\phi_n) \ge \dashint_{\M}  F(J\phi_n)  = & \frac{V_{A_n}}{V_{\M}} \dashint_{A_n}  F(J\phi_n)+\frac{V_{B_n}}{V_{\M}} \dashint_{B_n}  F(J\phi_n) =\frac{V_{A_n}}{V_{\M}} F\brk{\dashint_{A_n}  J\phi_n}+\frac{V_{B_n}}{V_{\M}} \dashint_{B_n}  F(J\phi_n)   \\
\stackrel{(1)}{\ge} & \frac{V_{A_n}}{V_{\M}} F\brk{\dashint_{A_n}  J\phi_n}+\frac{V_{B_n}}{V_{\M}} F\brk{\dashint_{B_n}  J\phi_n} \stackrel{(2)}{\ge} F\brk{\frac{V_{A_n}}{V_{\M}} \dashint_{A_n}  J\phi_n+\frac{V_{B_n}}{V_{\M}} \dashint_{B_n}  J\phi_n }\\ =&F\brk{\dashint_{\M}  J\phi_n  } = F\brk{\frac{V_{\phi_n\brk{\M} } }{V_{\M}}  }  \ge  F\brk{\frac{V_{\N}}{V_{\M}}}.
\end{split}
\eeq
where in the second equality on the first line we used the affinity of $F|_{[0,\frac{1}{4}]}$, and inequalities $(1),(2)$ are Jensen. In the last line, we used the injectivity of $\phi_n$ and the assumption $V_{\N} \le V_{\M}$.
\begin{lemma}  
\label{lem:asymptotic_Jacob_avg}
Suppose that  $V_{\N} \le V_{\M}$. 
Let $\phi_n \in \Lippi(\M,\N)$ and assume that $F\brk{\dashint_{\M}  J\phi_n} \to F\brk{ \frac{V_{\N}}{V_{\M}}}$. Then $ \dashint_{\M}  J\phi_n \to \frac{V_{\N}}{V_{\M}}$. Thus if we assume in addition that $\qrt < \frac{V_{\N}}{V_{\M}}$, then there exists $r_0 >\qrt$ such that $\dashint_{\M}  J\phi_n \ge r_0$ for sufficiently large $n$. 
\end{lemma}

\begin{proof} 
Since $\lim_{x \to \infty} F(x)=\infty$, $\dashint_{\M}  J\phi_n$ is bounded and we may assume it converges to $L$. $\dashint_{\M}  J\phi_n=\frac{V_{ \phi_{n}\brk{\M} }}{V_{\M}} \le \frac{V_{\N}}{V_{\M}}$ impies $ L \le \frac{V_{\N}}{V_{\M}}$. 
The assumption $F(\dashint_{\M}  J\phi_n) \to F\brk{ \frac{V_{\N}}{V_{\M}}}$ and the continuity of $F$ imply that $F(L)=F\brk{ \frac{V_{\N}}{V_{\M}}}$. Since $F|_{[0,1]}$ is strictly decreasing and  $0\le L \le \frac{V_{\N}}{V_{\M}} \le 1$, it follows that $L = \frac{V_{\N}}{V_{\M}} $. 
\end{proof}

\subsubsection{The Jacobian of a minimizing sequence is greater than $\bqrt$}
In the next proposition, we  do not assume anything on the areas of $\M,\N$. 
\begin{proposition}
\label{prop:concentrated_Jacobian}
Let $\phi_n \in \Lipp(\M,\N)$. 
Assume that $\dashint_{\M}  F(J\phi_n)-F\brk{\dashint_{\M}  J\phi_n  } \to 0$, $\dashint_{\M}  J\phi_n \ge r_0$ for some $r_0 >\qrt$, and that $\sup_{n}\dashint_{\M}  J\phi_n<\infty$. Then $\lim_{n \to \infty} V_{\{J\phi_n \le \frac{1}{4}\}}=0$. 
\end{proposition}
\begin{comment}
The assumption $\sup_{n}\dashint_{\M}  J\phi_n<\infty$ can be omitted, but we won't need this stronger version.
An assumption that ensures $\dashint_{\M}  J\phi_n >\bqrt$ is certainly needed; otherwise if e.g.~$J\phi_n < \bqrt$ everywhere on $\M$, then $\dashint_{\M}  F(J\phi_n)=F\brk{\dashint_{\M}  J\phi_n  }$ since $F|_{[0,\qrt]}$ is affine.
\end{comment}

\begin{proof}
Set again $A_n:=\{ p \in \M \, | \, J\phi_n(p) \le \frac{1}{4} \}$.
Setting 
\[
\lambda_n= \frac{V_{A_n}}{V_{\M}}, a_n=\dashint_{A_n}  J\phi_n,b_n=\dashint_{B_n}  J\phi_n,
\]
Inequality $(2)$ from Equation \eqref{eq:Euclideanquant1} is 
\[
 \lam_n F(a_n)+(1-\lam_n)F(b_n) \ge F\brk{\lam_n a_n+\brk{1-\lam_n}b_n}=F(c_n),
\] 
where 
\[c_n:=\lam_n a_n+(1-\lam_n)b_n=\dashint_{\M}  J\phi_n \ge r_0 > \qrt.
\]
The assumption  $\dashint_{\M}  F(J\phi_n)-F\brk{\dashint_{\M}  J\phi_n  } \to 0$ implies that 
\beq
\label{eq:Euclideanquant2}
D_n:= \lambda_nF(a_n)+(1-\lambda_n)F(b_n)-F\brk{\lambda_n a_n +\brk{1-\lambda_n}b_n} \to 0.
\eeq
We first show that $b_n$ is bounded. If not, then by passing to subsequences we may assume that $a_n \to a, c_n \to c, b_n \to \infty$. (By definition $a_n \le \bqrt$ is bounded, and $c_n$ is bounded by our assumption.)
Since 
\[
\lam_n a_n+(1-\lam_n)b_n=c_n=\dashint_{\M}  J\phi_n 
\]
is bounded, we must have $\lam_n \to 1$.  In fact we have
\beq
\label{eq:just_eq}
1-\lambda_n=\frac{c_n-a_n}{b_n-a_n}\sim (c-a)b_n^{-1}.
\eeq
Using $F(b_n)\geqslant F(c_n)+(b_n-c_n)F'(c_n)$ we get
\[
\begin{split}
&D_n+F(c_n)=\\
&\lambda_n F(a_n)+(1-\lambda_n)F(b_n) \ge \\
&  \lambda_n F(a_n)+(1-\lambda_n)F(c_n)+(1-\lambda_n)(b_n-c_n)F'(c_n). 
\end{split}
\]
Taking limits of both sides we obtain
\beq
\label{eq:just_eq2}
\liminf_{n \to \infty} D_n + F(c)\ge  F(a)+(c-a)F'(c).
\eeq
where the evaluation of the limit of the RHS follows from estimate \eqref{eq:just_eq}.
The strict convexity of $F|_{[\qrt,\infty)}$ implies 
\[
F(a)+(c-a)F'(c) > F(c).
\] Indeed, this is equivalent to 
$ \frac{F(c)-F(a)}{c-a}< F'(c)$. The mean value theorem implies that the LHS equals $F'(t)$ for some $t \in (a,c)$. Since $  c>\bqrt$, the strict convexity implies the strict inequality $F'(t) < F'(c)$.  So, \eqref{eq:just_eq2} implies $\liminf_{n \to \infty} D_n >0$ which contradicts \eqref{eq:Euclideanquant2}. 
We showed that $b_n$ is bounded, so we may assume that $a_n \to a, b_n \to b, \lam_n \to \lam$. 
$D_n \to 0$ then implies 
\[
\lam F(a)+(1-\lam)F(b)=F\brk{\lam a+\brk{1-\lam}b}.
\] Since $a \le \bqrt, b \ge c > \bqrt$ and $F|_{[\qrt,\infty)}$ is strictly convex we must have $\lam=0$ or $\lam=1$. $\lam=1$ is impossible, since it would imply $\lim_{n \to \infty}c_n=a \le \bqrt$, contradicting the assumption $c_n=\dashint_{\M}  J\phi_n \ge r_0 > \bqrt$.

\emph{A detailed proof that $\lam=0$:}
Let $\alpha\in (0,1)$ be such that 
\[
\alpha \frac{1}{4} + (1-\alpha) b = \lam a+(1-\lam)b.
\]
Then $\alpha \ge \lambda$ and from convexity we have
\[
F(\lam a+(1-\lam)b)=\lambda F(a) + (1-\lambda) F(b) \ge \alpha F(\frac{1}{4}) + (1-\alpha) F(b) \ge F(\alpha \frac{1}{4} + (1-\alpha) b ),
\]
which implies $\alpha F(\frac{1}{4}) + (1-\alpha) F(b) = F(\alpha \frac{1}{4} + (1-\alpha) b )$. The strict convexity of $F|_{[\qrt,\infty)}$ implies that  $\al=0$ or $\al=1$. Again $\al=1$ is excluded, since this would imply $\lim_{n \to \infty}c_n= \bqrt$. So $\al=0$, hence $\lam=0$ as well.
\end{proof}
\subsubsection{The norm of a minimizing sequence is concentrated at $J\phi_n	> \bqrt$}
Here we also do not assume anything on the areas of $\M,\N$.  
\begin{lemma}
\label{lem:concentration_ophi_norm}
Let $\phi_n \in \Lipp(\M,\N)$. Suppose that $E_2(\phi_n)- \dashint_{\M}  F(J\phi_n) \to 0$ and that $V_{A_n} \to 0$, where  $A_n=\{ J\phi_n \le \bqrt \}$. Then $\lim_{n \to \infty}   \int_{A_n}  |d\phi_n|^2=0$.
\end{lemma}
\begin{comment} 
We don't assume that $\dashint_{\M}  F(J\phi_n)-F\brk{\dashint_{\M}  J\phi_n  } \to 0$, nor that $\dashint_{\M}  J\phi_n  \ge r_0>\bqrt$. Since we are not in the context of \propref{prop:concentrated_Jacobian}, we explicitly assumed that $V_{A_n} \to 0$. 
\end{comment}

In order prove \lemref{lem:concentration_ophi_norm} we shall need the following pointwise estimate:
\begin{lemma}
 \label{lem:quantbound_vol_distortion_dist_SO22}
Let $A \in M_2$ with $\det A \ge 0$. Then
\beq
\label{eq:quantbound_vol_distortion_dist_SO2_1}
\dist^2(A,K) \le \dist^2(A,\SOtwo)-(1-2\det A) \le 2\dist^2(A,K).
\eeq
This is a quantitative generalization of \lemref{lem:bound_vol_distortion_dist_SO22} in the regime where $\det A \le \bqrt$.
The proof is given in \secref{proofest1}.
\end{lemma}

\begin{proof}[Of \lemref{lem:concentration_ophi_norm}]
\[
 \int_{\M}  \dist^2(d\phi_n,\SOtwo) -F(J\phi_n) =\int_{A_n}  \dist^2(d\phi_n,\SOtwo) -F(J\phi_n)+\int_{B_n}  \dist^2(d\phi_n,\SOtwo) -F(J\phi_n) \to 0.
 \]
Since both summands are nonnegative $\int_{A_n}  \dist^2(d\phi_n,\SOtwo) -(1-2J\phi_n) \to 0$. By \lemref{lem:quantbound_vol_distortion_dist_SO22} $\int_{A_n}  \dist^2(d\phi_n,K) \to 0.$ 

We note that for any $A \in M_2$, $|A| \le \dist(A,K)+1$. Indeed, $K$ is compact and $X \in K \Rightarrow |X| \le 1$, so if $X \in K$ satisfies $|A-X|=\dist(A,K)$, then
\[
|A| \le |A-X|+|X| \le \dist(A,K)+1.
\] 
Using $ (x+y)^2 \le 2(x^2+y^2),$ we deduce that $|A|^2 \le 2\dist^2(A,K)+2$. Thus
 \[
\int_{A_n}  |d\phi_n|^2=\int_{A_n}  \brk{|d\phi_n|^2-2} + 2V_{A_n} \le 2\int_{A_n}  \dist^2(d\phi_n,K) + 2V_{A_n} \to 0.
 \]  
  \end{proof}
 \subsubsection{A minimizing sequence is asymptotically conformal}

We prove that a minimizing sequence is asymptotically conformal. We use the notation $\CO:=\{ \lambda Q \, | \, \lam \ge 0, Q \in \SOtwo \}$ for the set of (weakly orientation-preserving) conformal matrices. We use the well-known fact that for $A \in M_2$ with $\det A \ge 0$ and singular values $\sigma_1 \le \sigma_2$, $\dist^2(A,\CO)=\frac{1}{2} ( \sig_1-\sig_2)^2.$ We shall need the following estimate (which we prove in Appendix \ref{conformal_estimate}):
\begin{lemma}
 \label{lem:quantbound_vol_distortion_dist_confSO22}
 Let $A \in M_2$ with $\det A \ge \bqrt$, and let $\sigma_1 \le \sigma_2$ be its singular values. Then 
\beq
\label{eq:quantbound_vol_distortion_dist_SO2_1}
\brk{\sqrt \sig_2-\sqrt \sig_1}^4
  \le \dist^2(A,\SOtwo)-2(\sqrt{\det A}-1)^2 \le 2\dist^2(A,\CO).
\eeq
The RHS is valid under the weaker assumption $\det A \ge 0$. This statement is a quantitative generalization of \lemref{lem:bound_vol_distortion_dist_SO22} in the regime where $\det A \ge \bqrt$. 
\end{lemma}
\begin{comment}
There is no $c>0$ satisfying 
\[
\dist^2(A,\SOtwo)-2(\sqrt{\det A}-1)^2 \ge c\dist^2(A,\CO).
\]
\end{comment}
The following proposition do not assume anything on the areas of $\M,\N$.    
\begin{proposition}
\label{prop:asymptotic_conformal}
Let $\phi_n \in \Lipp(\M,\N)$ be bounded in $W^{1,2}(\M,\N)$. Suppose that $E_2(\phi_n)- \dashint_{\M}  F(J\phi_n) \to 0$ and that $V_{\{ J\phi_n \le \frac{1}{4} \}} \to 0$. Then $\int_{\M}  \dist^2(d\phi_n,\CO) \to 0$. 
\end{proposition} 

\begin{proof}
First, since
  $2\dist^2(A,\CO)=|A|^2-2\det A \le |A|^2$,
 we deduce from \lemref{lem:concentration_ophi_norm} that
\beq
\label{eq:asymptot_quan1a}
\lim_{n \to \infty} \int_{A_n}  \dist^2(d\phi_n,\CO)=0.
\eeq
Consider 
\[
\int_{B_n}  \dist^2(d\phi_n,\SOtwo) -F(J\phi_n)=\int_{B_n}  \dist^2(d\phi_n,\SOtwo)- 2(\sqrt{J\phi_n}-1)^2 \to 0.
\]
\lemref{lem:quantbound_vol_distortion_dist_confSO22} implies that
\beq
\label{eq:middle_eq}
\int_{B_n} \brk{\sqrt{\sig_2(d\phi_n)}-\sqrt{\sig_1(d\phi_n)}}^4 \to 0.
\eeq
\[
\begin{split}
\int_{B_n}  \dist^2(d\phi_n,\CO) &=\half \int_{B_n}  \brk{\sig_1 -\sig_2}^2=\half \int_{B_n}  \brk{\sqrt{\sig_1} -\sqrt{\sig_2}}^2 \brk{\sqrt{\sig_1} +\sqrt{\sig_2}}^2  \\
& \le \half \sqrt{\int_{B_n} \brk{\sqrt{\sig_1} -\sqrt{\sig_2}}^4} \sqrt{\int_{B_n} \brk{\sqrt{\sig_1} +\sqrt{\sig_2}}^4}
\end{split}
\]
where we used Holder's inequality in the last step.
Since $\sqrt{\sig_i} \le \sqrt{|d\phi_n|} \Rightarrow \sqrt{\sig_1}+\sqrt{\sig_2} \le 2\sqrt{|d\phi_n|}$, we deduce that
\[
\int_{B_n}  \dist^2(d\phi_n,\CO)\le 2 \sqrt{\int_{B_n} \brk{\sqrt{\sig_1} -\sqrt{\sig_2}}^4} \sqrt{\int_{B_n} |d\phi_n|^2} \le  2 \|d\phi_n\|_{L^2}  \sqrt{\int_{B_n} \brk{\sqrt{\sig_1} -\sqrt{\sig_2}}^4}
\]
Since $\|d\phi_n\|_{L^2} $ is bounded, we conclude from Equation \eqref{eq:middle_eq} that 
\[
\lim_{n \to \infty} \int_{B_n}  \dist^2(d\phi_n,\CO)=0.
\]
Together with \eqref{eq:asymptot_quan1a} this implies the required assertion.
\end{proof}

We prove that the weak limit of an asymptotically conformal sequence is conformal.

In the next lemma, we assume $\M,\N \subseteq \mathbb{R}^2$.
\begin{lemma}
\label{lem:limi_asymptot_conf_is_conf}
Let $\M,\N \subseteq \mathbb{R}^2$ be Euclidean domains. Let $\phi_n \weakly{} \phi$ in $W^{1,2}(\M,\N)$, and suppose that $J\phi_n \ge 0$ a.e.~and $\dist^2(d\phi_n,\CO) \to 0$. 
Then $\phi$ is weakly-conformal, i.e.~$d\phi \in \CO$ a.e.~and $d\phi_n \to d\phi$ strongly in $L^2(K)$ for every $K \Subset \M^\circ$. 
\end{lemma}

It is not clear whether the conclusion $d\phi \in \CO$ a.e.~holds for maps between manifolds.
\begin{proof}
Let $K \Subset \M^\circ$. 
$\phi_n \weakly{} \phi$ in $W^{1,2}(\M,\N)$ implies that $ J\phi_n \rightharpoonup J\phi $ in $L^1(K)$. Using the fact that for any $A \in M_2$ with $\det A \ge 0$, $|A|^2-2\det A=2\dist^2(A,\CO)$, we get
\beq
\label{eq:some33}
\lim_{n\to \infty} \|d\phi_n\|_{L^2(K)}^2=2\lim_{n\to \infty} \int_K J\phi_n=2 \int_K J\phi \stackrel{(1)}{\le}\|d\phi\|_{L^2(K)}^2, 
\eeq
where in inequality $(1)$ we used $|A|^2\ge 2\det A$ for $A=d\phi$. 

Since $d\phi_n \rightharpoonup d\phi$ in $L^2$, and the $L^2$-norm is weakly lower semicontinuous, we deduce that $\|d\phi\|_{L^2(K)}=\lim_{n\to \infty} \|d\phi_n\|_{L^2(K)}$. The weak convergence+convergence of norms imply that $d\phi_n \to d\phi$ strongly in $L^2(K)$.


In particular, we have equality in inequality $(1)$ which implies that $\phi$ is weakly conformal. 
\end{proof}

We conclude this subsection by proving that the limit of asymptotically surjective maps is surjective.
\begin{lemma}
\label{lem:injeclimitinject}
Let $\phi_n \in \Lippi(\M,\N)$ satisfy $V_{\phi_n(\M)} \to V_{\N}$, and suppose that $\phi_n$ converges strongly in $W^{1,2}$ to a continuous function $\phi: \M \to \N$. Then $\phi$ is surjective. If $\phi \in \Lipp(\M,\N)$, it is injective a.e.
\end{lemma}


\begin{proof}
Assume by contradiction that $\exists q \in \N\setminus \phi(\M)$. Since $\phi(\M)$ is compact, it is closed in $\N$. $q \in \N \setminus \phi(\M)$ which is open, thus $\phi(\M)\subseteq \N \setminus U$ for some open neighbourhood $U\subset \N$ of $q.$ Choose a bump function $\alpha:\N \to \R$ with $\int_{\N} \alpha=1$ whose support lies in $U.$

The Jacobians $J\phi_n$ converge in $L^1$ to $J\phi.$  (This \emph{strong} convergence of the Jacobians holds also for maps between manifolds, see \cite{374383}.)

The maps $\alpha\circ \phi_n$ are uniformly bounded and converge in measure to $\alpha\circ \phi.$ Since product of a bounded sequence which converges in measure and a sequence that converges in $L^1$ converges in $L^1$, $(\alpha \circ \phi_n) J\phi_n \to (\alpha \circ \phi) J\phi$ in $L^1$. In particular,
\[
\int_{\M} (\alpha \circ \phi_n) J\phi_n\to \int_{\M} (\alpha \circ \phi) J\phi
\]
The RHS is zero because $\alpha\circ \phi\equiv 0.$ (since $\supp \al \subseteq U$ and $\phi(\M) \subseteq \N \setminus U$). But by the area formula, we have $\int_{\M} (\alpha \circ \phi_n) J\phi_n=\int_{\phi_n(\M)} \al \to \int_{\N} \alpha=1$, which is a contradiction.

Assume that $\phi \in \Lipp(\M,\N)$. Then since $\phi$ is surjective, 
\[
V_{\N}=\lim_{n \to \infty} V_{\phi_n(\M)}= \lim_{n \to \infty} \int_{\M} J\phi_n=\int_{\M} J\phi  =   \int_{\N} |\phi^{-1}(y)| \ge  \int_{\N} 1 =V_{\N}
\]
which implies that $|\phi^{-1}(y)|=1$ a.e.
\end{proof}

\subsubsection{Proof of \thmref{thm:asymptot_hom}}


The assumption $ E_2(\phi_n) \to  F\brk{ \frac{V_{\N}}{V_{\M}}}$  implies that $\phi_n$ is bounded in $W^{1,2}(\M,\N)$. (by Poincare's inequality). Thus we may assume that it converges weakly in $W^{1,2}$ to some $\phi$. Equation \eqref{eq:Euclideanquant1} implies that $\dashint_{\M}  F(J\phi_n)-F\brk{\dashint_{\M}  J\phi_n  } \to 0$ and that $E_2(\phi_n)- \dashint_{\M}  F(J\phi_n) \to 0$. Since $ E_2(\phi_n) \ge F\brk{\dashint_{\M}  J\phi_n}$ and $\lim_{x \to \infty} F(x)=\infty$, $\dashint_{\M}  J\phi_n$ is bounded from above.
By \lemref{lem:asymptotic_Jacob_avg} $\dashint_{\M}  J\phi_n \ge r_0 > \bqrt$. Thus, by \propref{prop:concentrated_Jacobian}
 $ V_{\{J\phi_n \le \frac{1}{4}\}} \to 0$, so \propref{prop:asymptotic_conformal} implies $\dist^2(d\phi_n,\CO) \to 0$. 

\lemref{lem:limi_asymptot_conf_is_conf} then implies that $\phi$ is weakly conformal and $d\phi_n \to d\phi$ strongly in $L^2(K)$ for every $K \Subset \M^\circ$. 
The strong convergence $d\phi_n \stackrel{L^2(K)}{\to} d\phi$ implies 
 \[
 E_2(\phi|_{K}) = \lim_{n \to \infty}E_2(\phi_n|_{K}) \le\lim_{n \to \infty}E_2(\phi_n)=F\brk{ \frac{V_{\N}}{V_{\M}}} .
 \]
 So, we proved that $ E_2(\phi|_{K})  \le F\brk{ \frac{V_{\N}}{V_{\M}}}$ for every $K \Subset \M^\circ$. By the monotone-convergence theorem, 
 \beq
 \label{eq:ineq_fin_energy}
  E_2(\phi)  \le F\brk{ \frac{V_{\N}}{V_{\M}}}. 
  \eeq
 Furthermore, since $J\phi_n \ge 0$ and $ J\phi_n \weakly{L^1(K)} J\phi $ , 
 \[
  \int_{K}  J\phi= \lim_{n \to \infty} \int_{K}  J\phi_n \le \lim_{n \to \infty} \int_{\M}  J\phi_n=V_{\N},
  \] where the last equality is by \lemref{lem:asymptotic_Jacob_avg}.
$J\phi_n \ge 0$ and  $ J\phi_n \weakly{L^1(K)} J\phi $ imply $J\phi \ge 0$ a.e.~By monotone convergence $ \int_{K}  J\phi \le V_{\N}$  for every $K  \Subset \M^\circ$ implies $ \int_{\M}  J\phi \le V_{\N}$, or $ \dashint_{\M}  J\phi \le \br \le 1$. Thus
\beq
\label{eq:fin_fin}
  F\brk{ \frac{V_{\N}}{V_{\M}}}  \ge  E_2(\phi)   \ge \dashint_{\M} F(J\phi) \ge F\brk{ \dashint_{\M} J\phi} \ge   F\brk{ \frac{V_{\N}}{V_{\M}}},
\eeq 
where we used Inequality \eqref{eq:ineq_fin_energy} and the monotonicity of $F|_{[0,1]}$. Thus $E_2(\phi)=F\brk{ \frac{V_{\N}}{V_{\M}}} $, and $\dashint_{\M} J\phi=\frac{V_{\N}}{V_{\M}} > \qrt$. We are now in the same situation as in the proof of  \thmref{the:main_Euclid_bound}; \lemref{lem:Jensen_equality1} implies that $J\phi$ is constant a.e., hence $J\phi=\frac{V_{\N}}{V_{\M}} $. We already proved that $\phi$ is conformal, so it is a homothety. (The conformality of $\phi$ follows also from the equality $E_2(\phi)  = \dashint_{\M} F(J\phi)$.) 

In particular we have 
\[
\lim_{n \to \infty} \dashint_{\M}  J\phi_n=\frac{V_{\N}}{V_{\M}} = \dashint_{\M} J\phi.
\] Repeating the same argument below Equation \eqref{eq:some33} with $K$ replaced by $\M$, we deduce that $d\phi_n \to d\phi$ strongly in $L^2(\M)$, and hence by Poincare's inequality, $\phi_n \to \phi$ strongly in $W^{1,2}(\M,\N)$.
\lemref{lem:injeclimitinject} then implies that $\phi$ is surjective and injective almost everywhere. (For deducing the almost injectivity we are using the fact that $\phi$ is Lipschitz; $\phi|_{\M^\circ}$, being a homothety, is smooth by known regularity results, so it is Lipschitz on $\M^\circ$ hence extends to a Lipschitz map on $\M$. Alternatively, we can also argue that $\phi \in W^{1,\infty}$, hence Lipschitz).

It follows that $\phi$ is injective on $\M^\circ$.
Indeed, assume $\phi(p_1) =q = \phi(p_2)$, where $p_1 \neq p_2\in \M^\circ$ and $q\in \N^\circ$. (since $d\phi$ is invertible $\phi(\M^\circ) \subseteq \N^\circ$, see e.g.~Ex 4.2 in \cite{Lee13}.)
Since $\phi|_{\M^\circ}:\M^\circ \to \N^\circ$ is a local diffeomorphism, there exist disjoint open neighborhoods $U_i\ni p_i$ and $V\ni q$ such that $\phi(U_i) = V$, hence $\Vol\brk{ \{ q\in \N : |\phi^{-1}(q)|>1 \}} \ge \Vol(V) > 0, $
which is a contradiction to $\phi$ being injective a.e.~(We proved here that every injective a.e.~smooth map from a manifold without boundary into a manifold, having invertible differential is injective). This completes the proof of the first part of the theorem.

\emph{Proof that under additional assumptions $\phi$ is a diffeomorphism}

Since $\phi_n\to \phi$ in $W^{1,2}(\M;\N)$, it follows that $\phi_n|_{\pl\M}\to \phi|_{\pl\M}$ in $L^2(\pl\M;\pl\N)$, and (after taking a subsequence) pointwise almost everywhere in $\pl\M$. Since $\phi_n(\pl \M)\subset \pl \N$, and since $\pl\N$ is closed and $\phi$ is continuous we conclude that $\phi(\pl\M)\subset \pl \N$. We already established that $\phi(\M)=\N$, $\phi(\M^\circ) \subseteq \N^\circ$, which together with $\phi(\pl\M)\subset \pl \N$ imply that $\phi(\M^\circ) = \N^\circ,\phi(\pl\M)=\pl \N$.
In particular, $\phi|_{\M^\circ}:\M^\circ \to \N^\circ$ is a diffeomorphic homothety.
Thus for every  $x,y\in\M$, let $\M^\circ\ni x_n\to x$ and $\M^\circ\ni y_n\to y$; then
\[
d_{\N}(\phi(x),\phi(y)) =
\lim_{n \to \infty} d_{\N}(\phi(x_n),\phi(y_n)) =
\lim_{n \to \infty} \lam d_{\M}(x_n,y_n) =\lam d_{\M}(x,y).
\]
This implies that $\phi$ is injective on all of $\M$. Finally, by a version of the Myers-Steenrod theorem for manifolds with boundary (for a short argument see \cite{253994}), if $\pl \M,\pl \N$ are smooth, then $\phi$ is smooth up to the boundary. 
\subsection{Proof of \thmref{thm:low_regime_rigidity} }

Given $\phi \in \Lipp(\M,\N)$,\lemref{lem:quantbound_vol_distortion_dist_SO22} implies that 
\[
\dist^2(d\phi,K) \le \dist^2(d\phi,\SOtwo)-(1-2J\phi) \le 2\dist^2(d\phi,K).
\]
Integrating we get
\beq
\label{eq:quantbound_vol_distortion_dist_SO2_1a}
\dashint_{\M} \dist^2(d\phi,K)  \le E_2(\phi) - \brk{ 1-2 \dashint_{\M} J\phi} \le  2\dashint_{\M} \dist^2(d\phi,K).
\eeq
Applying this for $\phi \in \Lippi(\M,\N)$, under the assumption $V_{\phi\brk{\M} }/ V_{\M}  \le \bqrt$ proves assertion \eqref{eq:quantbound_vol_distortion_dist_SO2_1abc}. Suppose now that $\br \le \bqrt$; then 
\[
\begin{split}
E_2(\phi) - F\brk{\frac{V_{\N}}{V_{\M}}}& =E_2(\phi) -  F\brk{\frac{V_{\phi\brk{\M} }}{V_{\M}}} +  F\brk{\frac{V_{\phi\brk{\M} }}{V_{\M}}}-F\brk{\frac{V_{\N}}{V_{\M}}}\\
&=E_2(\phi) -  F\brk{\frac{V_{\phi\brk{\M} }}{V_{\M}}} + \frac{2(V_{\N}-V_{\phi(\M)}) }{V_{\M}}. 
\end{split}
\]
where we used the fact that $F(x)-F(y) =2(y-x)$ for $x,y \in [0,\qrt]$.
This proves assertion \eqref{eq:quantbound_vol_distortion_dist_SO2_1abcde}. 

\section{Geometric properties of the minimizing well}
\subsection{Geometric characterization of the well $K$}
\label{sec:secgeom_prop_well}
We give an alternative description of the well $K$ defined in \ref{def:doublewell_combined}.
Let $\GLtwo$ be the group of real $2 \times 2$ matrices having positive determinant. 
\begin{proposition}
\label{prop:geometric_cof_char_well}
Let $A \in \GLtwo$.  The following are equivalent:
\begin{enumerate}
\item $A \in K$. ($\sig_1(A)+\sig_2(A)=1$.) 
\item $\Cof A + A \in \SOtwo$.
\item $\Cof A + A =O(A)$,
where $O(A)$ is the orthogonal polar factor of $A$, i.e.~$A=O(A)P(A)$, where $O(A) \in \SOtwo$ and $P(A)$ is symmetric positive-definite. We note that $O(A)$ is also the closest matrix to $A$ in $\SOtwo$. 
\end{enumerate}
$\Cof A$ here denotes the standard cofactor matrix of $A$. 

Moreover, if $A -O(A)=\alpha \Cof A$ for some $\al \in \R$, then either $\al=-1$ (and then $A \in K$) or $A$ is conformal with singular values $\sig_1=\sig_2=\frac{1}{1-\al}$.
\end{proposition}

\begin{proof}
All three conditions are invariant under left and right multiplication by special orthogonal matrices; this follows from the multiplicative properties $\Cof (AB)=\Cof A \Cof B$, $O(UAV)=UO(A)V$ for any $U,V \in \SOtwo$, together with the fact that $\Cof Q =Q$ for every $Q \in \SOtwo$.
Using SVD, this reduces the problem to showing equivalence of the conditions for the special case where $A=\Sigma $ is diagonal positive-definite. 
In that case
\[
\Sigma+\Cof \Sigma =\begin{pmatrix} \sig_1 & 0 \\\ 0 & \sig_2 \end{pmatrix}+\begin{pmatrix} \sig_2 & 0 \\\ 0 & \sig_1 \end{pmatrix}=(\sig_1+\sig_2)\id,
\]
so 
\[
\Sigma+\Cof \Sigma \in \SOtwo \iff \Sigma+\Cof \Sigma=\id=O(\Sigma) \iff \sig_1+\sig_2=1.
\]

Now, assume that $A -O(A)=\alpha \Cof A$ for some $\al \in \R$.  Again, the orthogonal invariance of the equation implies that  $\Sigma -O(\Sigma)=\alpha \Cof \Sigma$ or $\Sigma -\id=\alpha \Cof \Sigma$, where $\Sigma$ is the diagonal part in the SVD of $A$. Thus
\[
\sig_1-1=\al \sig_2, \sig_2-1=\al \sig_1.
\]
By subtracting we deduce that $(\al+1)(\sig_1-\sig_2)=0$, so  either $\al=-1$ or $\sig_1=\sig_2$. 
\end{proof}

\subsection{Maps in the well $K$ are critical}

We prove \propref{prop:wellK_is_critic} which states that maps in the well $K$ are critical points of the energy. Let $(\M,\g),(\N,\h)$ be smooth $n$-dimensional Riemannian manifolds, $\M$ compact. Let $\phi \in C^2(\M,\N)$ with $J\phi >0$.
Let $O(d\phi):T\M \to \phi^*T\N$ be the (unique) closest orientation-preserving \emph{isometric} section to $d\phi$. (This is essentially the orthogonal polar factor of $d\phi$).

As we show in Appendix \ref{sec:EL_deriv}, The Euler-Lagrange equation of the functional 
\[
E_2(\phi)=\int_{\M} \dist^2\brk{d\phi,\SO(\g,\phi^*\h)}=\int_{\M}|d\phi-O(d\phi)|^2
\]  is
\beq
\label{eq:EL_elastic}
\delta \brk{d\phi-O(d\phi)}=0,
\eeq
and the one for $E_p$ is
\beq
\label{eq:EL_elastic_gen_p}
\delta \brk{\dist^{p-2}\brk{d\phi,\SO(\g,\phi^*\h)} \brk{d\phi-O(d\phi)}}=0,
\eeq
where the coderivative $\delta_{\nabla^E}: \Omega^1(\M;\phi^*T\N) \to  \Gamma(\phi^*T\N)$
is the adjoint of the connection $\nabla^{\phi^*T\N}:\Gamma(\phi^*T\N)\to \Omega^1(\M;\phi^*T\N)$. (For precise definitions see e.g.~\cite{EL83}. In the Euclidean case, where $\M=\N=\R^n$ are endowed with the standard metrics, $\delta$ is the standard row-by-row divergence.)  

We recall the following statement (the Piola identity):
\begin{quote}
\emph{
For every $C^2$ map $\phi:\M \to \N$,  $\, \, \, \delta \brk{\Cof d\phi }=0$,
}
\end{quote}
where the cofactor of $d\phi$ is defined intrinsically using the metrics on $\M,\N$; see \cite{KS19} (Section 2.1) for details. (If $\M=\N=\R^n$, $\Cof d\phi$ is the standard cofactor matrix.) The Piola identity, which holds for any sufficiently regular map, is well-known in the Euclidean case. (see e.g.~\cite[Ch.~8.1.4.b]{Eva98} and \cite[p.~39]{Cia88} for a proof.) It was generalized to mappings between arbitrary Riemannian manifolds in \cite{KS19}.

Considering Equation \eqref{eq:EL_elastic_gen_p}, the Piola identity suggests a natural way to find critical points of $E_p$---look for maps $\phi$ which satisfy 
\beq
\label{eq_cof_spec}
\dist^{p-2}\brk{d\phi,\SO(\g,\phi^*\h)} \brk{d\phi-O(d\phi)}=\al \Cof d\phi,
\eeq
where $\al \in \R$ is constant. We specialize to the case where $\M,\N$ are $2$-dimensional.

\begin{proof}[of \propref{prop:wellK_is_critic}]

By \propref{prop:geometric_cof_char_well} $d\phi \in K \Rightarrow d\phi-O(d\phi)=- \Cof d\phi$, hence the Piola identity implies that $\phi$ satisfies \eqref{eq:EL_elastic}.
If the singular values of $d\phi$ are constant, then $\dist\brk{d\phi,\SO(\g,\phi^*\h)}$ is constant, so $\phi$ satisfies \eqref{eq:EL_elastic_gen_p} as well.

On the other hand, let $p \neq 2$ and suppose that $d\phi \in K$ and that $\phi$ satisfies \eqref{eq:EL_elastic_gen_p}. Then 
\beq
\label{eq:midd_EL_const}
0=\delta \brk{H\brk{d\phi-O(d\phi)}}=\delta \brk{H \Cof d\phi}=\delta \brk{\Cof d\phi},
\eeq
where $H:=\dist^{p-2}\brk{d\phi,\SO(\g,\phi^*\h)}$.
The last equality implies that $H$ is constant: Write $\omega:=\Cof d\phi \in \Omega^1(\M;\phi^*T\N)$. Given an orthonormal frame $E_i$ for $T\M$, we have
\beq
\delta \omega  =-\trg(\nabla \omega)=-\sum_{i=1}^d (\nabla_{E_i} \omega)(E_i),
\label{eq:codiv_formula1}
\eeq
where $\nabla \omega$ is the connection induced  on $T^*\M \otimes \phi^*T\N$ by the Levi-Civita connection on $\M$ and $\nabla^{\phi^*T\N}$ (see e.g.~\cite[Lemma 1.20]{EL83} for a proof). 

Equations \eqref{eq:codiv_formula1} and \eqref{eq:midd_EL_const} imply that
\[
0= \delta(H \omega )= H\delta(\omega) - \tr_{\g}(dH \otimes \omega)=-\tr_{\g}(dH \otimes \omega),
 \]
 so
 \[
 0=\tr_{\g}(dH \otimes \omega)=\sum_{i=1}^d dH(E_i)\cdot \omega(E_i)=\sum_{i=1}^d dH(E_i)\cdot \Cof d\phi(E_i).
 \]
Since $\Cof d\phi:T\M \to \phi^*T\N$ is invertible, $dH=0$, so $H$ is constant.

Thus 
\[
 \sig_1(d\phi)+ \sig_2(d\phi)=1, \,\,\,  \brk{\sig_1(d\phi)-1}^2+\brk{\sig_2(d\phi)-1}^2
\]
are constant, hence $ \sig_1(d\phi), \sig_2(d\phi)$ are constant.
\end{proof}

The following result states that the suggested approach based on the Piola identity in \eqref{eq_cof_spec} produces only homotheties or maps in the well $K$.
\begin{lemma}
\label{lem:Piola_approach_limit}
Let $\phi \in C^2(\M,\N)$ with $J\phi>0$. Then $d\phi-O(d\phi)=\al \Cof d\phi$ if and only if $d\phi \in K$ or $\phi$ is a homothety with $\sig(d\phi)=\frac{1}{1-\al}$.

For $p \neq 2$, $\dist^{p-2}\brk{d\phi,\SO(\g,\phi^*\h)} \brk{d\phi-O(d\phi)}=\al \Cof d\phi$ if and only if $\phi$ has constant singular values.
\end{lemma}

The proof is given in in Appendix \ref{sec:add_proof}.


Another application of \propref{prop:wellK_is_critic} is the following corollary: 
We say a map $\phi \in C^2(\M,\N)$ is \emph{affine} if $\nabla d\phi = 0$, where $\nabla =\nabla^{T^*\M \otimes \phi^*T\N}$ is the natural connection induced on $T^*\M \otimes \phi^*T\N$ by the Levi-Civita connections on $\M,\N$. This notion of affinity coincides with the standard one in the Euclidean case.
\begin{corollary}
\label{lem:non_affine_crit}
Let $\M,\N$ be any Riemannian surfaces. Then there exist local maps $\M \to \N$, which are critical points of $E_p$, for every $p \ge 1$. 
In particular, there exist non-affine critical points of $E_p$, for every $p \ge 1$. 
\end{corollary}

By local existence, we mean that given any two points $p \in \M,q \in \N$ there exist open neighbourhoods $U,V$ of $p,q$ respectively, and a map $\phi:U \to V$ sending $p$ to $q$, which is $E_p$-critical.

\begin{proof}
Given positive numbers $\sig_1<\sig_2$ there always exist locally many maps $\M \to \N$ with constant singular values $\sig_1,\sig_2$, see discussion of the Euclidean case in \cite{351550} and for the general Riemannian case in \cite{383251}. If we take singular values satisfying $\sig_1+\sig_2=1$, then \propref{prop:wellK_is_critic} implies that such maps are critical.
\end{proof}

\paragraph{Non-regular solutions}
The flexibility of the well $K$ implies that the solutions to the Euler-Lagrange equation \eqref{eq:EL_elastic} do not have to be regular. 

Let $\M=\N=\R^2$ endowed with the standard Euclidean metrics.
%
%
Since $K_{\sig_1,\sig_2}$ is rank-one connected for $\sig_1<\sig_2$ (see \cite[p. 190]{desimone2002macroscopic}), there exist non-differentiable maps $\phi \in W_{loc}^{1,\infty}(\R^2,\R^2)$ which satisfy $d\phi \in K_{\sig_1,\sig_2}$ a.e.~

If $\sig_1+\sig_2=1$, then by \propref{prop:geometric_cof_char_well} $d\phi-O(d\phi)=- \Cof d\phi$ a.e..
A weak version of the Piola identity implies that any such map satisfies a weak form of Equation \eqref{eq:EL_elastic_gen_p}.  Thus, weak solutions of \eqref{eq:EL_elastic_gen_p} are not necessarily $C^1$. 

\subsection{Critical maps having constant singular values are in $K$}
We prove \propref{prop:constant_sing}, which states that critical maps having constant singular values between Euclidean spaces, are either affine or in the well $K$.
Note that a map having constant singular values that is $E_p$ critical for some $p$, is $E_q$ critical for any value of $q$, and in particular $E_2$-critical.

\begin{proof}
Let $\phi \in C^2(\Omega,\R^2)$ have constant singular values. If $\sig_1=\sig_2$, then $\phi$ is a homothety and in particular affine, so there is nothing to prove. Suppose that $\sig_1 \neq \sig_2$.
Let $d\phi=U\Sigma V^T$ be the SVD of $d\phi$, $U,V \in \SOtwo$ .Since $\sig_1 \neq \sig_2$, $U,V$ can be chosen smoothly, locally around every point $x \in \Omega$, see e.g. \cite{3163368}.


Let $\diag_2$ be the vector space of real-valued two-by-two matrices. 
Consider the map $T: \diag_2 \to C^{\infty}(\Omega,\R^2)$ given by
\[
T:\begin{pmatrix} a & 0 \\\ 0 & b \end{pmatrix}  \to
\div \brk{U\begin{pmatrix} a & 0 \\\ 0 & b \end{pmatrix} V^T},
\]
where $\div$ is the divergence operator, acting row-by-row.
We shall use the following result (which we prove below):
\begin{lemma}
\label{lem:constancy}
If $T=0$, then $U,V$ are constant.
\end{lemma}

The Euclidean Piola identity is 
\[
\div (\Cof d\phi)=0.
\]
Since $\Cof d\phi=U\Cof \Sigma V^T$, this translates into $T(\Cof \Sigma)=0$.
Now assume that $\phi$ is $E_2$ critical, i.e.~ 
\[
\div \brk{d\phi-O(d\phi)}=0.
\]
Since $d\phi-O(d\phi)=U(\Sigma-\id)V^T$, this translates into $T(\Sigma-\id)=0$. We established $\Cof \Sigma, \Sigma-\id \in \ker T$. Since $\diag_2$ is two-dimensional, if $\Cof \Sigma, \Sigma-\id$ are linearly independent, then $T=0$, hence by \lemref{lem:constancy} $U,V$ are constant, and $\phi$ is affine.

If $\Cof \Sigma, \Sigma-\id$ are dependent, \propref{prop:geometric_cof_char_well} implies that either $d\phi \in K$ or $\phi$ is a homothety (and in particular affine).
\end{proof}

\begin{proof}[of \lemref{lem:constancy}]
Write $U=\begin{pmatrix} c & -s \\\ s & c \end{pmatrix}, V^T=\begin{pmatrix} \til c & - \til s \\\ \til s & \til c \end{pmatrix}$.
Writing explicitly 
\[
U\begin{pmatrix} a & 0 \\\ 0 & b \end{pmatrix} V^T= \begin{pmatrix} ac\til c-b s\til s & -ac\til s-b\til c s \\\ a s\til c+b \til s c & -a s \til s+b c \til c \end{pmatrix},
\]
so
\[
\begin{split}
T\brk{\begin{pmatrix} a & 0 \\\ 0 & b \end{pmatrix}}&=\begin{pmatrix}  a \bdx(c\til c)-b \bdx(s\til s) -a\bdy(c\til s)-b\bdy(\til c s)   \\\ a \bdx(s\til c)+b\bdx( \til s c) -a \bdy(s \til s)+b \bdy(c \til c) \end{pmatrix} \\
&=\begin{pmatrix}  a \brk{\bdx(c\til c) -\bdy(c\til s)} +b\brk{ -\bdx(s\til s) -\bdy(\til c s) }  \\\ a \brk{ \bdx(s\til c)- \bdy(s \til s)}  +b \brk{ \bdx( \til s c) + \bdy(c \til c)} \end{pmatrix}\\
&= \begin{pmatrix} \bdx(c\til c) -\bdy(c\til s) & -\bdx(s\til s) -\bdy(\til c s) \\\ \bdx(s\til c)- \bdy(s \til s) & \bdx( \til s c) + \bdy(c \til c) \end{pmatrix} \begin{pmatrix} a   \\\ b \end{pmatrix}
\end{split}
\]
Thus $T=0$ implies
\[
 \begin{pmatrix} \bdx(c\til c) -\bdy(c\til s) & -\bdx(s\til s) -\bdy(\til c s) \\\ \bdx(s\til c)- \bdy(s \til s) & \bdx( \til s c) + \bdy(c \til c) \end{pmatrix}=0.
\]
Rewriting this we get the following
\beq
\label{eq:cs_1}
\begin{split}
&\partial_x(\tilde sc)=-\partial_y(c\tilde c) \\
 &\partial_y(\tilde s c)=\partial_x(c\tilde c), \\
\end{split}
\eeq
\beq
\label{eq:cs_2}
\begin{split}
&\partial_x(s\tilde c)=\partial_y(s\tilde s) \\
&\partial_y(s\tilde c)=-\partial_x(s\tilde s).
\end{split}
\eeq

The system \eqref{eq:cs_1} implies that $f=-\tilde sc+i(c\tilde c)=i\cos(\theta)e^{i\tilde \theta}$ is holomorphic, and the system \eqref{eq:cs_2} implies that $g=s\tilde c+i(s\tilde s)=\sin(\theta)e^{i\tilde \theta}$ is holomorphic. 
Thus $(if)^2+g^2=e^{2i\tilde \theta}$ is holomorphic. Since its image lies on a circle, the open mapping theorem implies that it is constant, thus $e^{i\tilde \theta}$ is constant. 

Together with the holomorphicity of $f,g$, we deduce that $c=\cos(\theta), s=\sin(\theta)$ are holomorphic, so are also constant.
\end{proof}
\section{Energy minimizers of other functionals}
\label{sec:gen_distortion_functional_prf}
%
Throughout the following subsection, we assume the setting described in \secref{sec:Other_distortion_functionals}. In particular, $E=E_f$ and $F$ are defined as in equations  \eqref{eq:gen_energy_functional} and \eqref{eq:min_uniform0}.

We prove in \lemref{lem:minfunc_prop} in \secref{sec:add_proof} that $F$ is well-defined, and also show it has the same monotonicity properties as $f$ does. Our starting point in the analysis is the observation that the convexity property of $F$ is the key element in passing from the pointwise bound to a variational bound. This element was also present in \secref{sec:Euclidean_distortion_functional}. 

We shall need the following localized notions of convexity:

\begin{definition}
Let $f:(a,b) \to \R$, and let $c \in (a,b)$. We say that $f$ is convex at $c$ if $f\left(\alpha x + (1- \alpha)y \right) \leq \alpha f(x) + (1-\alpha)f(y)$
holds whenever $ \alpha \in [0,1]$ and $x,y \in (a,b)$ satisfy $\alpha x + (1- \alpha)y =c$. Strict pointwise convexity is defined similarly.

Furthermore, we say that $f$ is midpoint-convex at $c$ if $f(\frac{x + y}{2}) \le \frac{f(x) + f(y)}{2}  $
holds whenever $x,y \in (a,b)$ satisfy $\frac{x + y}{2} =c$.

We note that midpoint-convexity plus continuity at a point does not imply convexity at a point, even though midpoint-convexity plus continuity on an \emph{interval} does imply (full) convexity at that interval.

\end{definition}

Our basic observation is the following: 
\begin{lemma}
\label{lem:convex_F_implies_varia_bound}
Let $\M,\N$ be compact Riemannian surfaces, and let $\phi \in \Lipp(\M,\N)$. Suppose that $F$ is convex at $\dashint_{\M} J\phi$. Then $E(\phi) \ge F\brk{\dashint_{\M} J\phi}.$
If $F$ is strictly convex at $\dashint_{\M} J\phi$, then equality implies that $J\phi$ is constant. 
\end{lemma}

\begin{proof}
The definition of $F$ together with the convexity assumption implies that  
\beq
E(\phi) \stackrel{(1)}{\ge} \dashint_{\M} F(J\phi) \stackrel{(2)}{\ge} F\brk{\dashint_{\M} J\phi}.
\eeq
\end{proof}

When $F$ is affine on a subinterval the Jacobian of energy minimizers need not be constant (we saw that already in the special case of the Euclidean functional in \thmref{the:main_Euclid_bound_injective}).  This observation inspires the problem of characterizing the cost functions $f$ which give rise to $F$ which have affine parts. This seems a non-trivial problem, and the only such case we are aware of is when $f(x)=(x-1)^2$. Similarly, it is not clear which $f$'s give rise to convex $F$ .
Indeed, replacing the quadratic penalty with cubic or quartic penalties $f(x)=|x-1|^3, f(x)=(x-1)^4$ makes $F$ non-convex (see \cite{opt20}).
\subsection{Analysis of when $(\sqrt s,\sqrt s)$ is a minimizer}
\label{sec:minpoint1}
We are looking for conditions on $f$ which ensure that the solution to the minimization problem \eqref{eq:min_uniform0} is obtained at $(\sqrt s,\sqrt s)$, that is
\beq
\label{eq:min_uniform1}
F(s)=\min_{xy=s,x,y>0} f(x)+ f(y)=2f(\sqrt s),
\eeq
As mentioned in the proof of $\lemref{lem:minfunc_prop}$, the minimum point for $s \in (0,1)$ is obtained when $x,y \le 1$. Thus, Equation \eqref{eq:min_uniform1} is equivalent to 
\beq
\label{eq:min_uniform2}
f(\sqrt{xy}) \le \frac{f(x) + f(y)}{2} \, \,\, \, \text{ for every } \, \,\, x,y \in (0,1], \, \, xy=s.
\eeq
After defining $g:\R \to [0,\infty)$, by $g(x) = f(e^x)$, or $f(x) = g(\log x)$,
\eqref{eq:min_uniform2}  becomes
\[
 g\left( \frac{\log x + \log 
y}{2}\right) \le  \frac{g(\log x) + g(\log y)}{2} \, \,\, \, \text{ whenever } \, \, \log x, \log y \in (-\infty,0], \, \,  \log x+ \log y=\log s.
\]
Equivalently,  $g|_{(-\infty,0]}$ is midpoint-convex at $\half \log s$, i.e.
\[
g\brk{\frac{x + y}{2}} \le \frac{g(x) + g(y)}{2}  \, \,\, \, \text{ whenever } \, \,  x,  y \in(-\infty,0], \, \,   x+  y=\log s.
\]
Thus, we proved the following

\begin{lemma}
\label{lem:uniq_min_single_ratio}
Let $s \in (0 ,1]$.  $(\sqrt s,\sqrt s)$ is a (unique) minimizer of \eqref{eq:min_uniform0} if and only if $g|_{(-\infty,0]}$ is (strictly) midpoint-convex at $\half \log s$.
\end{lemma}

Note that $g|_{(-\infty,0]}$ is midpoint-convex at $\half \log s$ if and only if $g|_{[\log s,0]}$ is midpoint-convex at $\half \log s$.

\begin{corollary}
\label{cor:uniq_min_range}
Let $\de \in (0 ,1]$ and consider the following statements:
\begin{enumerate}
\item $g|_{[\log \de,0]}$ is (strictly) convex.
\item $(\sqrt s,\sqrt s)$ is a (unique) minimizer of \eqref{eq:min_uniform0} for every $s \in [\de ,1]$.
\item $g|_{[\half \log \de,0]}$ is (strictly) convex. 
\end{enumerate}
Then $(1) \Rightarrow (2) \Rightarrow (3)$. The converse implications do not hold in general.

In particular, $(\sqrt s,\sqrt s)$ is a (unique) minimizer of \eqref{eq:min_uniform0} every $s \in (0,1]$ if and only if $g|_{(-\infty,0]}$ is (strictly) convex. 
\end{corollary} 

\begin{proof}
$(1) \Rightarrow (2)$ follows directly from \lemref{lem:uniq_min_single_ratio}. 

$(2) \Rightarrow (3)$: By \lemref{lem:uniq_min_single_ratio},  $g|_{(-\infty,0]}$ is (strictly) midpoint-convex at every $\half \log s \in [\half \log \de,0]$. In particular, $g|_{[\half \log \de ,0]}$ is (strictly) midpoint-convex, and (strict) midpoint-convexity plus continuity implies full (strict) convexity. We do not provide examples that refute the converse implications.
\end{proof}




Since convex functions whose derivative obtains negative values tend to $+\infty$ when $x \to -\infty$, we get the following:

\begin{corollary}
\label{cor:min_uni_diverge}
Suppose that $(\sqrt s,\sqrt s)$ is a minimizer of \eqref{eq:min_uniform0} for every $s \in (0,1)$. Then $\lim_{x\to 0} f(x)=\infty$. 
\end{corollary} 
\corref{cor:min_uni_diverge} can be strengthened as follows:
\begin{proposition}
 \label{prop:div_inf_min}
 Suppose that $f$ does not diverge to infinity at zero. Then there exists $\ep>0$ such that $(\sqrt s,\sqrt s)$ is not a minimizer of \eqref{eq:min_uniform0} for every  $s \in (0,\ep)$.
 \end{proposition}
 We prove \propref{prop:div_inf_min} in \secref{sec:add_proof}. We note that \corref{cor:uniq_min_range} implies the following: 

\begin{corollary}
\label{cor:uniq_min_close}
Suppose that $f$ is differentiable and not flat at $x=1$. There exists $\de \in (0,1)$ such that $(\sqrt s,\sqrt s)$ is the unique minimizer of \eqref{eq:min_uniform0} for $s \in (\de, 1]$.  
\end{corollary}

\begin{proof}

The relation $g(x)=f(e^x)$ implies that $f$ is not flat at $x=1$ if and only if $g$ is not flat at $x=0$. If $g$ is not flat at $x=0$, it is strictly convex in some neighbourhood of it, and we can apply \corref{cor:uniq_min_range}.
\end{proof}


\subsection{Convexity of $F$}

As we saw in \lemref{lem:convex_F_implies_varia_bound}, the convexity of $F$ plays a key role in the analysis.
An interesting observation is that $(\sqrt s,\sqrt s)$ being a minimizer has implications on convexity properties of $F$.  In the following let $f,g$ be as in \secref{sec:minpoint1}, i.e.~$f(x)=g(\log x)$.

\begin{lemma}
\label{lem:short_convex_2a}
Let $0<\de<1$. If $g|_{[\log \de,0]}$ is strictly convex, then $F|_{[\de,1]}$  is strictly convex.
\end{lemma}
\begin{proof}
By \corref{cor:uniq_min_range}, $F(s)=2f(\sqrt s)=2g(\half \log s)$ for every $ s \in [\de,1]$.   

The function $s \to g(\half \log s)$ is strictly convex on $[\de,1]$, since it is a composition of the strictly concave function $s \to \half \log s$ together with the strictly decreasing and strictly  convex function $g|_{[\half \log \de,0]}$. Thus, $F|_{[\de,1]}$ is strictly convex. 

\end{proof}

In particular, we obtain the following:

\begin{corollary}
\label{cor:min_uniform_conv}
Let $g:\R \to [0,\infty)$ be a continuous function, which is strictly decreasing and strictly convex on $(-\infty,0]$, and strictly increasing on $[0,\infty)$, with $g(0)=0$. Set $f(x) = g(\log x)$, and define $F: (0,\infty) \to  [0,\infty)$ as in \eqref{eq:min_uniform0}. 

Then for every $s \in (0,1)$, $(\sqrt s,\sqrt s)$ is the unique minimizer of \eqref{eq:min_uniform0}, and $F|_{(0,1]}$  is strictly convex. 

\end{corollary}
%

We shall use the following lemma, whose proof we postpone to Appendix (in \secref{sec:convexity_results}).
\begin{lemma}
\label{lem:short_convexity}
Let  $F:(0,\infty) \to [0,\infty)$ be a continuous function which is left-differentiable at $x=1$, satisfying $F(1)=F_-'(1)=0$, that is strictly increasing on $[1,\infty)$, and strictly decreasing on $(0,1]$. Suppose that $F|_{[1-\epsilon,1]}$ is (strictly) convex for some $\epsilon>0$. Then there exists $\delta>0$ such that $F$ is (strictly) convex \emph{at every point} $y \in (1-\delta,1]$.
 \end{lemma}
The assumption $F_-'(1)=0$ cannot be omitted. 
This lemma lifts the convexity on a subinterval to a global convexity- global in the sense that convexity at a point is a statement about far away points, which lie outside the interval where convexity is initially given.


\lemref{lem:short_convex_2a} and  \lemref{lem:short_convexity} imply the following

\begin{corollary}
\label{cor:convex2ba_loc_to_glob}
Suppose that $f$ is differentiable and not flat at $x=1$.

 There exists $\epsilon \in (0,1)$ such that $F$ is strictly convex at every point $y \in (1-\epsilon,1]$.
\end{corollary}
\begin{proof}
Since $g(x)=f(e^x)$ is not flat at $x=0$ and attains a minimum there, $g|_{[\log \de,0]}$ is strictly convex for some $0<\de<1$. \lemref{lem:short_convex_2a} implies that $F|_{[\de,1]}$  is strictly convex.

By \corref{cor:uniq_min_range}, $F(x)=2f(\sqrt x)$ on $[\de,1]$.
Since $f$ attains a minimum at $x=1$, we have $F_-'(1)=f'(1)=0$. By \lemref{lem:minfunc_prop} $F$ satisfies all the conditions in \lemref{lem:short_convexity}, which implies the assertion.





%
 \end{proof}


\subsection{Variational Bounds}

We combine all our preliminary results in order to prove the variational claims.
\begin{proof}[of \thmref{thm:bijective_homo_min_short}]

\corref{cor:convex2ba_loc_to_glob} implies that $F$ is strictly convex at every point $y 
\in [\al,1]$ for some $\al \in (0,1)$. 

If $ \dashint_{\M} J\phi \ge \al$, then  $F$ is strictly convex at  $\dashint_{\M} J\phi$, so by Jensen inequality

\beq
E(\phi) \stackrel{(1)}{\ge} \dashint_{\M} F(J\phi) \stackrel{(2)}{\ge} F\brk{\dashint_{\M} J\phi}  = F\brk{\frac{V_{\phi\brk{\M} }}{V_{\M}}}  \stackrel{(3)}{\ge} F\brk{\frac{V_{\N}}{V_{\M}}}.
\eeq

The strict convexity of $F$ at $\dashint_{\M} J\phi$ implies that inequality $(2)$ is an equality if and only if  $J\phi $ is constant. By \corref{cor:uniq_min_close} we can assume that $\al$ is sufficiently large such that $(\sqrt s,\sqrt s)$ is the unique minimizer of \eqref{eq:min_uniform0} for $s \in [\al,1]$. Since $J\phi \ge \al $ is constant, inequality $(1)$ is an equality  if and only if $\sigma_1(d\phi)=\sigma_2(d\phi) $ a.e.~Thus $ E(\phi) = F\brk{\frac{V_{\N}}{V_{\M}}}$ implies that $\phi$ is a homothety. Inequality $(3)$ is an equality if and only if $\phi$ is surjective. 

Suppose that $ \dashint_{\M} J\phi < \al$. Since $F$ is convex at $\al$, it has a supporting line at $\alpha$, i.e.
$F(x) \ge T(x):=m(x-\al)+F(\al)$ for some $m \in \mathbb{R}$, and for every $x>0$.

This implies
\[
E(\phi) \ge \dashint_{\M} F(J\phi) \ge \dashint_{\M} T(J\phi)=  T\brk{\dashint_{\M} J\phi} > T(\al) =F(\al)  \ge F\brk{\frac{V_{\N}}{V_{\M}}}, 
\]
where the last inequality is due to the assumption  $\al \le \frac{V_{\N}}{V_{\M}} \le 1$.
The strict inequality $T\brk{\dashint_{\M} J\phi} > T(\al)$ follows from $ \dashint_{\M} J\phi < \al$, together with the fact that the slope $m<0$. ($m<0$ since 
$T(1) \le F(1)=0<F(\al)=T(\al)$.)
%
\end{proof}

\begin{proof}[Of \thmref{thm:inject_gen_bound}]

\lemref{lem:minfunc_prop} and \corref{cor:min_uniform_conv} imply that $F$ satisfies the assumptions of \lemref{lem:Jens_variant}. Thus, 
\beq
E(\phi) \ge \dashint_{\M} F(J\phi) \ge F\brk{\dashint_{\M} J\phi} = F\brk{\frac{V_{\phi\brk{\M} }}{V_{\M}}} \ge F\brk{\frac{V_{\N}}{V_{\M}}}.
\eeq

where we have used the monotonicity of $F$ and $0 \le \frac{V_{\phi\brk{\M} }}{V_{\M}} \le \frac{V_{\N}}{V_{\M}} \le 1$. 

Now, assume $E(\phi) = F\brk{\frac{V_{\N}}{V_{\M}}}$. 
   Then $V_{\phi\brk{\M} }=V_{\N}$ and  $\dashint_{\M} F(J\phi) = F\brk{\dashint_{\M} J\phi}$. \lemref{lem:Jens_variant} implies that $J\phi $ is constant, hence $J\phi =\frac{V_{\N}}{V_{\M}}$ a.e.~on $\M$. Finally, 
   $f\brk{\sigma_1(d\phi)}+ f\brk{\sigma_2(d\phi)}= F(J\phi)$, so $\brk{\sigma_1(d\phi),\sigma_2(d\phi)} $ is a minimizer of \eqref{eq:min_uniform0} with $s=J\phi$. By \corref{cor:min_uniform_conv}, this implies $\sigma_1(d\phi)=\sigma_2(d\phi) $ so $\phi$ is a homothety.
\end{proof}

\subsubsection{A motivating example-logarithmic distortion}
\label{subsec:logarithmic_example}
We shall now consider the special case where $f(x)=(\log x)^2,g(x)=x^2$. As we explain below, this example has a geometric origin.

Let $I \in M_2$ be the identity matrix, and let $g_{I}$ be the standard Euclidean metric on the space of matrices, i.e.
\[
g_I(X,Y) = \tr(X^TY), \qquad  \text{where }  \, X,Y \in T_{I}\GLtwo \simeq M_2,
\]
Now, let $g$ be the Riemannian metric on $\GLtwo$ obtained by left-translating $g_I$, i.e.~$g$ is the unique left-invariant metric on $\GLtwo$ whose restriction to $T_{I}\GLtwo$
is $g_I$. $g$ induces a distance function $d_g$ (a metric in the sense of metric spaces) on $\GLtwo$, in the usual way --- the distance between any two points is the length of a minimizing geodesic between these points.

Given any distance function $d$ on $\GLtwo$, one gets a notion of ``distance from being an isometry" by setting
\[
 \distSO{A}=\inf_{Q \in \SOn} d(A,Q).
\]
Even though it is not known how to compute explicitly the distance $d_g(A,B)$ between arbitrary two points $A,B \in \GLtwo$, there is a formula for $\dist_g(A,\SOtwo)$ for an arbitrary matrix $A \in \GLtwo$. This problem was analyzed in a series of papers by Neff and co-workers \cite{neff2014logarithmic,neff2016geometry,lankeit2014minimization} (and a simplified proof was obtained in \cite{KS16}). 
The formula is given by 
\[
 \dist_g^2(A,\SOtwo)=\| \log \sqrt{A^TA}\|^2=\brk{\log \sig_1\brk{A}}^2+\brk{\log \sig_2\brk{A}}^2,
 \]
 where $\log \sqrt{A^TA}$ is the unique symmetric logarithm of the positive-definite matrix $A^TA$, and $\sig_i(A)$ are the singular values of $A$.
 
Thus, for a map $\phi:\M \to \N$ between Riemannian surfaces, we have
\[
 \dist_g^2(d\phi,\SOtwo)=\brk{\log \sig_1(d\phi)}^2+\brk{\log \sig_2(d\phi)}^2,
 \]
so
\[
E_f(\phi)= \int_{\M} f\brk{\sigma_1(d\phi)}+ f\brk{\sigma_2(d\phi)}=\int_{\M} \dist_g^2(d\phi,\SOtwo).
\]
for the choice of $f(x)=(\log x)^2$. This shows that this specific cost function arises naturally when one considers distortion functionals that are induced by Riemannian metrics with given symmetries.

In particular, we obtain the following corollary of \thmref{thm:inject_gen_bound}:

\begin{corollary}
Set
\[
E_f(\phi)= \int_{\M} \brk{\log \sig_1\brk{d\phi}}^2+\brk{\log \sig_2\brk{d\phi}}^2=\int_{\M} \dist_g^2(d\phi,\SOtwo).
\]
Then the homotheties, if they exist, are the unique energy minimizers among all injective maps.
\end{corollary}

\begin{comment}
In this case $F(s)=2g(\half \log s)=\half \log^2 s$ is not convex on the entire interval $(0,
\infty)$. Indeed, $F''<0$ on $(e,
\infty)$. This an example where $F|_{(0,1]}$ is convex but $F|_{(0,\infty)}$ is not convex. This is in contrast with the Euclidean case, where $F$ (from Definition \eqref{def:energy_point_vol_bd}) was convex on $(0,\infty)$. This shows the importance of our strong formulation of \lemref{lem:Jens_variant}, which only relies on the convexity of $F|_{(0,1]}$. (\lemref{lem:Jens_variant} is a key element in the proof of \thmref{thm:inject_gen_bound}).
\end{comment}

\section{Discussion}
\label{sec:Discuss}


\paragraph{Phase transitions for general functionals}

\thmref{thm:bijective_homo_min_short} states that for sufficiently close homothetic manifolds, the homotheties are the energy minimizers. 
\corref{cor:uniq_min_close} and \propref{prop:div_inf_min} explain why a phase-transition is expected when $f$ does not diverge at zero: when $s \to 1$, $(\sqrt s,\sqrt s)$ is the  minimizer of problem \eqref{eq:min_uniform0}, and when $s \to 0$ it stops being a minimizer.  
Lifting this ``pointwise" phase transition to the variational problem $\inf_{\phi \in \Lippi(\M,\N)} E(\phi)$ is left for future works. We stress again that $F$ may be non-convex, e.g.~for
$f(x)=|x-1|^3$, see \cite{opt20};  in such cases the convex envelope of $F$ should play a role in the analysis. 
\paragraph{Optimal compression of surfaces}
Fix a surface $\M$ and $\lam <1$. A natural problem is finding the optimal way to squeeze it into a surface $\N$ satisfying $V_{\N}=\lam^2 V_{\M}$, i.e.
\beq
\label{eq:gen_min_surf}
E_{\M,\lam}:=\inf_{\{ \N \, | \, V_{\N}=\lam^2 V_{\M}\}} E_{\M,\N}=\inf_{\{ \N \, | \, V_{\N}=\lam^2 V_{\M}\}} \inf_{\phi \in \Lippi(\M,\N)}  E_p(\phi) \ge F^{p/2}\brk{\lam}.
\eeq
By \thmref{the:main_Euclid_bound_injective}, for $\lam \ge 1/2$ there is a unique solution $\N=\lam \M$ and a unique optimal embedding (homothety).
%
For $\lam < 1/2$, the optimal target shapes $\N$ which attain $E_{\M,\lam}$ are not unique. Take e.g.~$\M=[-a,a]^2$, and let $\phi_c:\M \to \R^2$ be as in  \eqref{eq:logarithmic_min}; $\phi_c$ has singular values $\sig_1(c),\sig_2(c)$ satisfying $(\sig_1(c)+\sig_2(c))^2=4+c^2$; setting $\lam_c=\frac{1}{\sqrt{4+c^2}}$, the map
\[
\lam_c\phi_c: \M \to \N_c:=\lam_c\phi_c(\M)
\] 
has constant singular values which sum up to $1$. By \thmref{the:main_Euclid_bound_injective}, $E_p(\lam_c\phi_c) =  F^{p/2}\brk{\lam_c^2}$, hence 
\[
E_{\M,\lam_c}= \min_{\{ \N \, | \, V_{\N}=\lam_c^2 V_{\M}\}} E_{\M,\N}=E_{\M,\N_c}
\]
is attained at $\N=\lam_c\phi_c(\M)$; the square $[-a,a]^2$ is optimally mapped into a \emph{twisted} shrinked shape $\lam_c \phi_c([-a,a]^2)$ (Figure 1).
\begin{figure}[H]
  \captionsetup[subfigure]{labelformat=empty}
    \centering
    \subfloat[$c=1$]{{\includegraphics[width=3.5cm]{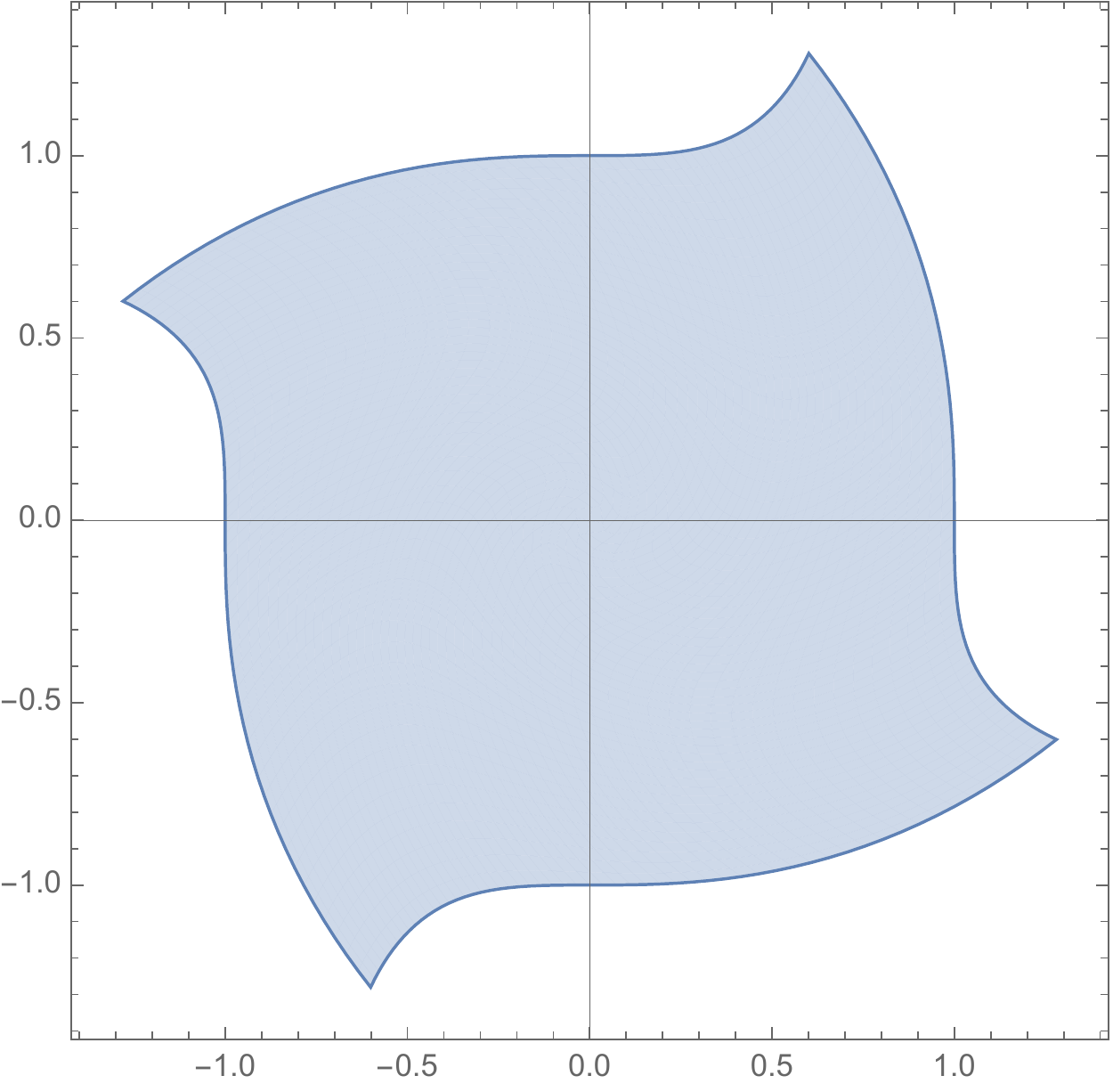} }}
    \qquad 
    \subfloat[$c=5$]{{\includegraphics[width=3.5cm]{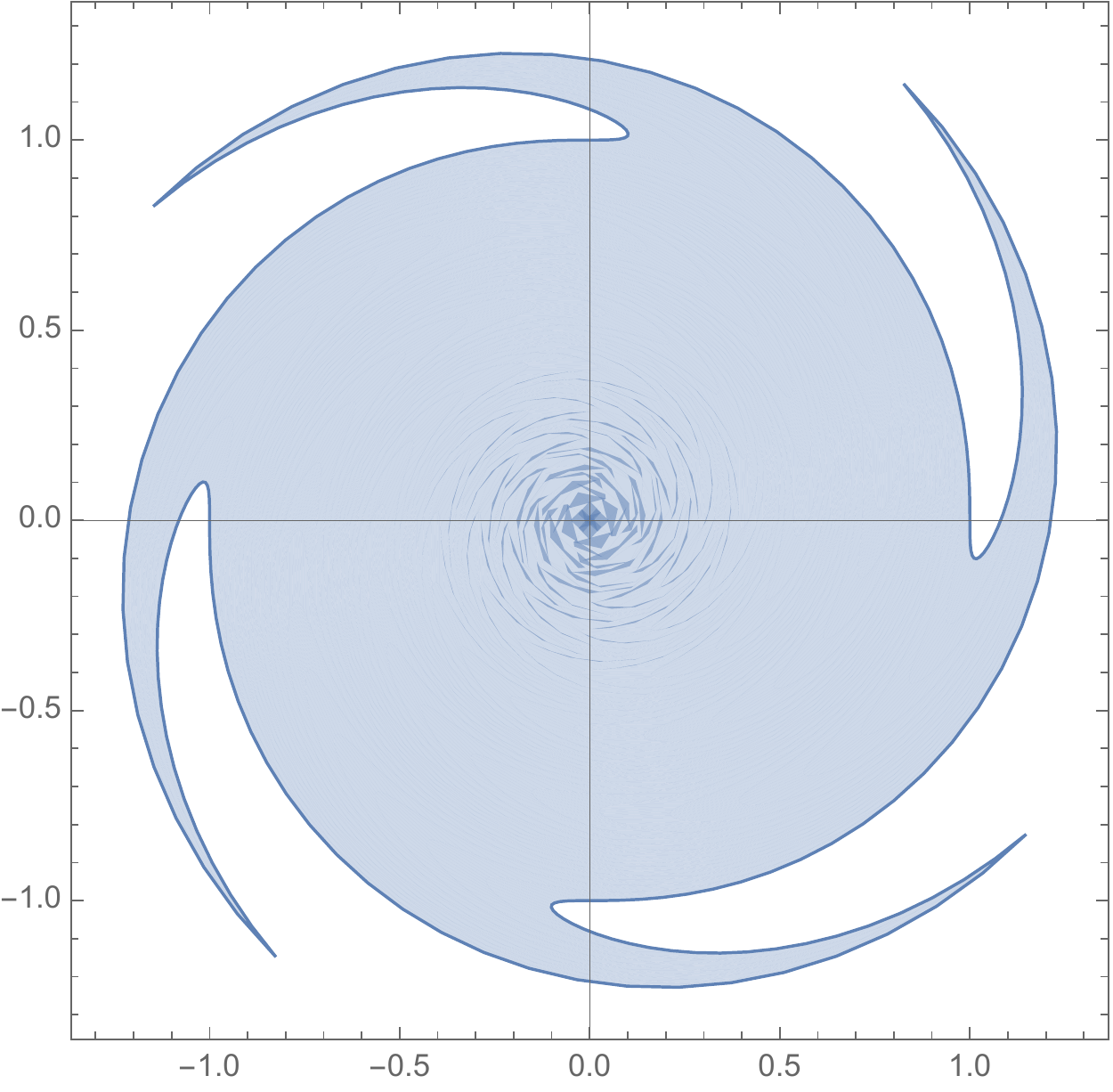} }}
     \qquad 
       \subfloat[$c=14$]{{\includegraphics[width=3.5cm]{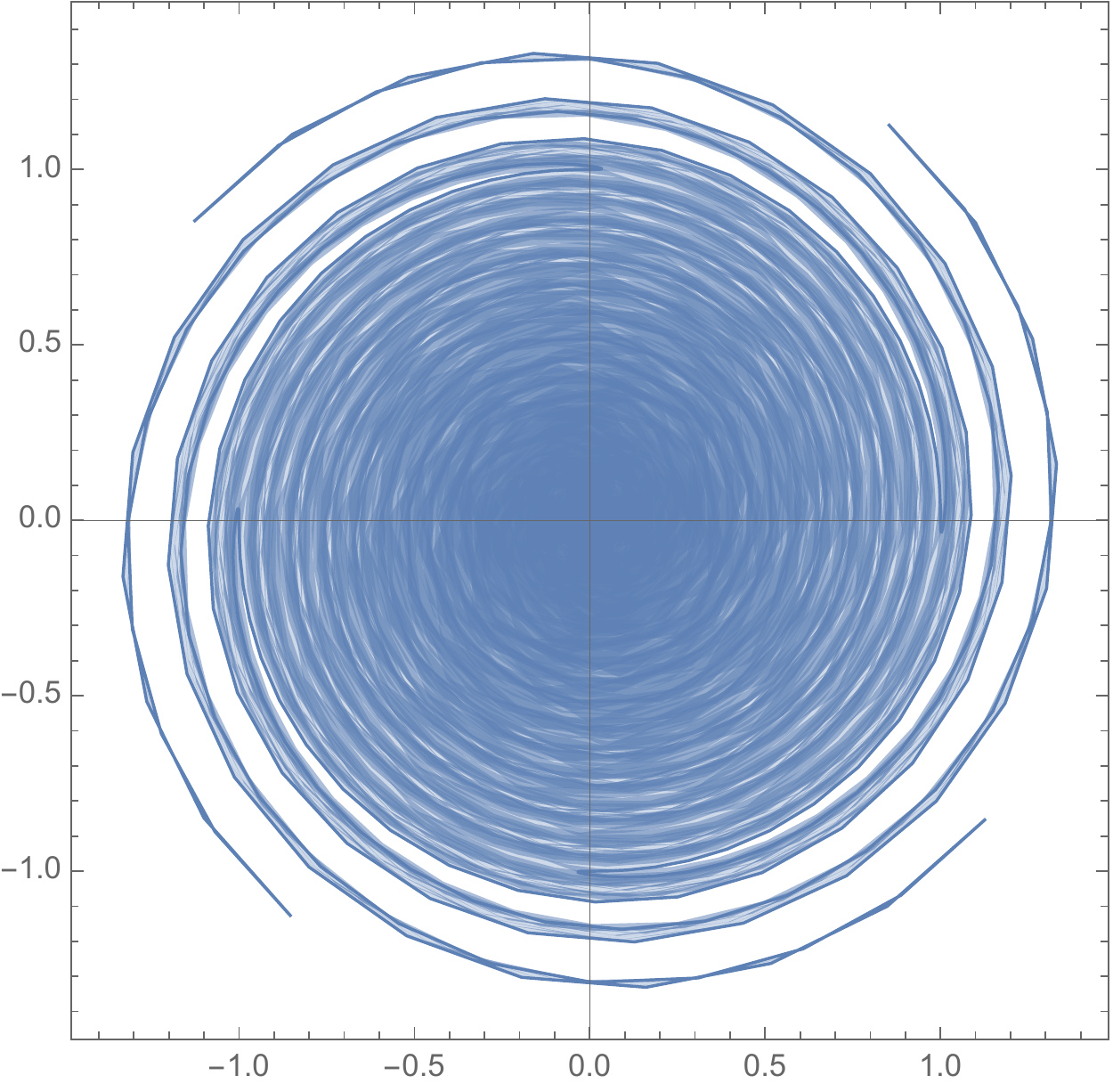} }} \\
              \subfloat[$c=51$]{{\includegraphics[width=3.5cm]{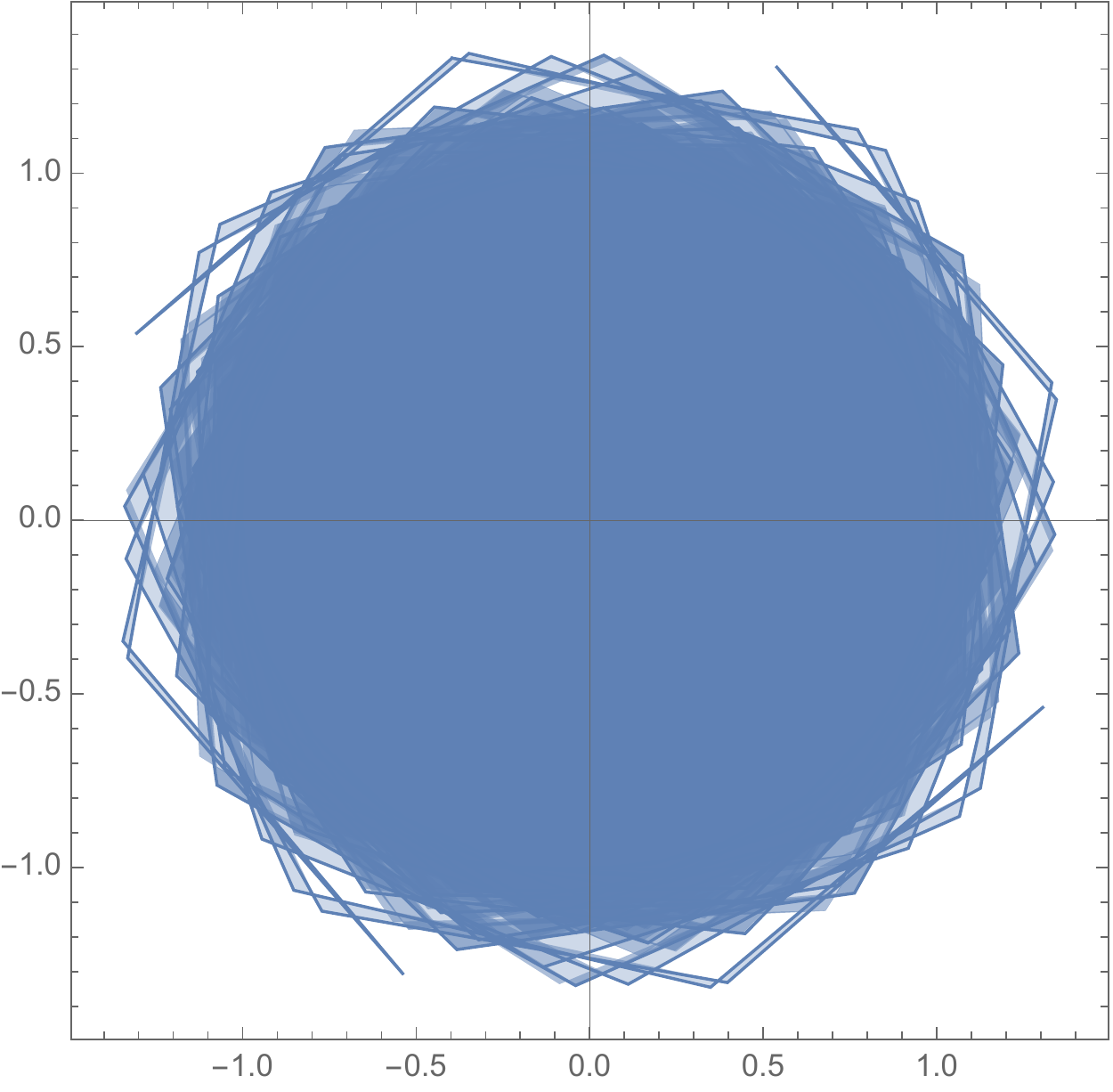} }}
                  \qquad 
                                \subfloat[$c=700$]{{\includegraphics[width=3.5cm]{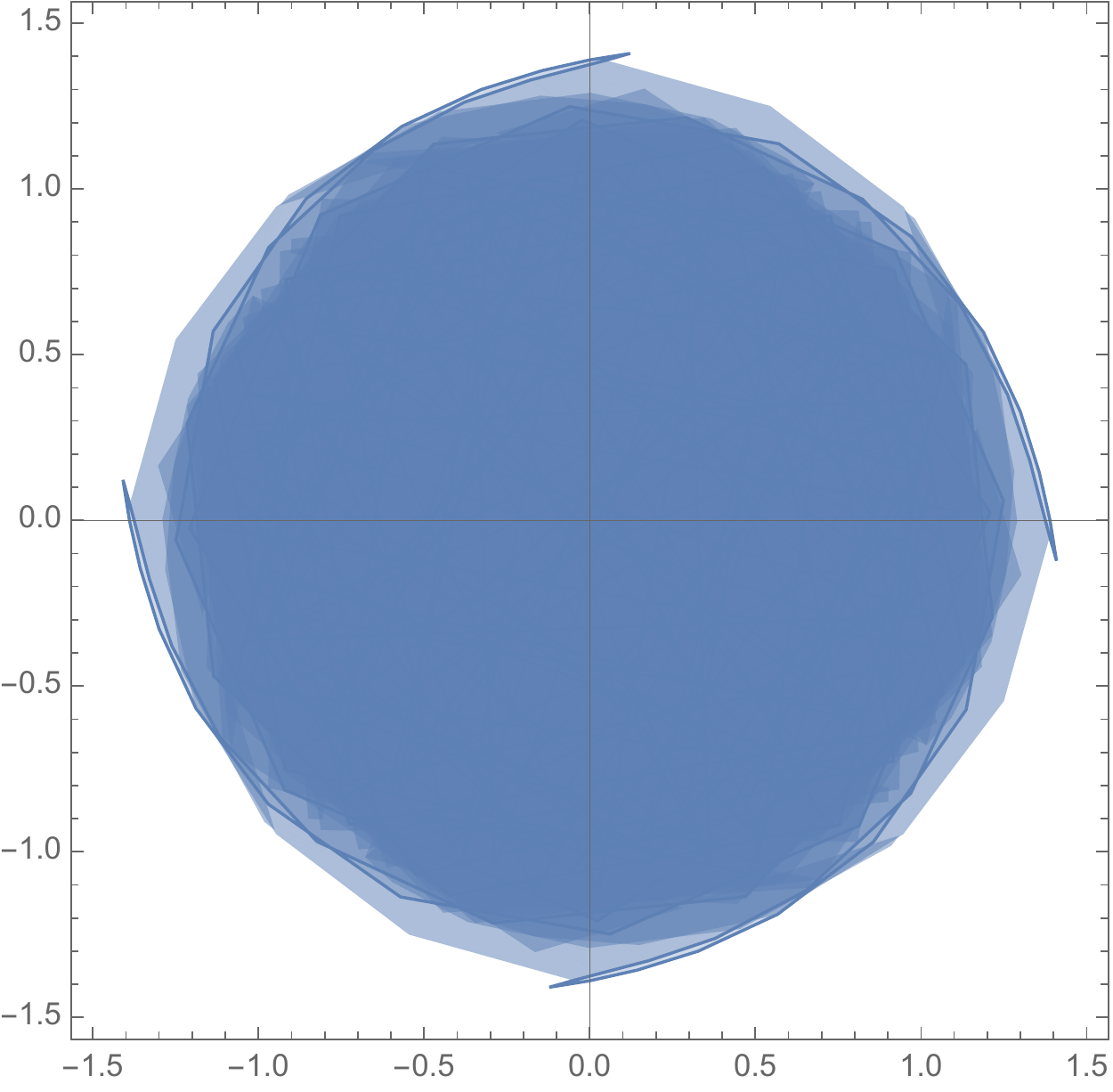} }}
    \caption{Shapes $\phi_c([-1,1]^2)$ for different values of $c$}
\end{figure}
%
%
A natural question is whether 
\[
E_{\M,\lam}= \min_{\{ \N \, | \, V_{\N}=\lam^2 V_{\M}\} } E_{\M,\N}= F^{p/2}(\lam^2)
\]
is attained at $\N=\lam \M$.
By \thmref{thm:equiv_min_self_map}, this is equivalent to the question whether there exists an area-preserving map $[-a,a]^2 \to [-a,a]^2$ having constant sum of singular values $1/\lam$. 

{\bfseries Acknowledgements}
We thank Stefan M\"uller for suggesting the proof that asymptotically	 conformal maps converge to a conformal map.
We thank Connor Mooney for suggesting the use of a concave function $\psi$ in \propref{prop:many_minimizers_peq2}, and Dmitri Panov for suggesting the area-preserving flow example $\phi_c$ in \eqref{eq:logarithmic_min}.  We thank Fedor Petrov for a suggested proof of a convexity result.
We thank Nadav Dym for suggesting a proof that no phase transition occurs for some energy functionals.
We thank Cy Maor for providing many helpful insights along the way. Finally, we thank Raz Kupferman for carefully reading this manuscript, and for suggesting various improvements during the research process.

This research was partially supported by the Israel Science Foundation (Grant No. 1035/17), and by a grant from the Ministry of Science, Technol- ogy and Space, Israel and the Russian Foundation for Basic Research, the Russian Federation.

\appendix


\section{Convexity results}
\label{sec:convexity_results}

%
%
\begin{proof}[of \lemref{lem:Jensen_equality1}]
The strict convexity of $F$ in $[a,\infty)$ implies that for $x>a\ge y$,
\[
F\brk{t x + (1-t) y} = tF(x) + (1-t)F(y) \Rightarrow t\in \{0,1\}.
\]
Denote $A = g^{-1}((a,\infty))$. If $\mu(A)=0$ or $\mu(A^c)=0$, we are done: If  $\mu(A)=1$, then $g>a$ a.e., so we are in the domain where $F$ is strictly convex. If $\mu(A)=0$ then $g \le a$ a.e.~- and the only case we need to check is when $\int_X g  \in [a,\infty)$. We then must have $\int_X g =a$, so $g=a$ a.e.

We show that $0<\mu(A)<1$ cannot occur. Denote $x=\dashint_A g$, and $y=\dashint_{A^c} g$; then $y\le a<x$ and $\int_X g= \mu(A)x + \mu(A^c)y$. Thus, using the equality assumption and Jensen's inequality,
\[
F(\int_X g) = \int_X F\circ g = \int_A F\circ g + \int_{A^c} F\circ g \ge \mu(A)F(x) + \mu(A^c)F(y)\ge F(\int_X g).
\]
Therefore equality holds, so the comment at the beginning of the proof implies that either $\mu(A)=1$ or $\mu(A^c)=1$. 

\end{proof}

The following lemma is a variant of Jensen inequality, when we have convexity only on a partial subset of our domain. 

\begin{lemma}
\label{lem:Jens_variant}
Let  $F:(0,\infty) \to [0,\infty)$. Assume that $F|_{(0,1]}$ is strictly decreasing and convex, with $F^{-1}(0)=\{1\}$. 


Let $g:X \to (0,\infty)$ be a measurable function defined on a probability space $X$ with  $\int_X g \in (0,1]$. Then $F(\int_X g) \le\int_X F\circ g\,\,$ and if equality occurs $g(x) \in (0,1]$ a.e.~ If $F|_{(0,1]}$ is strictly convex equality occurs if and only if $g(x)$ is constant a.e.

Equivalently: $F$ is convex at every point in $[0,1]$.

\end{lemma}

\begin{proof}
Set $F^*(x)=F(\min(x,1))$; $F^*\le F$ pointwise, and $F(x)=F^*(x)$ for $x\in[0,1]$.
Since  $\int_X g\in(0,1]$, $F\left(\int_X g\right)=F^*\left(\int_X g\right).$
Since $F^*$ is convex
\[
F\left(\int_X g\right)=F^*(\int_X g)\le \int_X F^*{\circ} g \le\int_X F{\circ} g
\]
as desired.
If there is an equality, then we have $F^*{\circ} g=F {\circ} g$ a.e., so $g \in (0,1]$ a.e.~If $F|_{(0,1]}$ is strictly convex then equality occurs if and only if $g(x)$ is constant a.e.~

We prove that $F^*$ is convex: (Just draw it!:)

We need to show that $F^*(t x + (1-t) y) \le tF^*(x) + (1-t)F^*(y)$.
If $t x + (1-t) y \ge 1$, then the LHS vanishes, so we are done. Thus, suppose that $t x + (1-t) y \le 1$. If $x,y$ are both not greater than $1$, then the assertion is just the convexity of $F|_{(0,1]}$. 

So, we may assume W.L.O.G that $0<x\le1< y$, and $t x + (1-t) y \le 1$. The assumptions 
imply $1\ge t x + (1-t) y \ge  t x + (1-t) 1 $, hence
\[
F\brk{t x + (1-t) y} \le F\brk{ t x + (1-t) 1} \le  tF(x) + (1-t)F(1) =tF(x),
\]
where the second inequality is due to the convexity of $F|_{(0,1]}$. Thus
\[
F^*\brk{ t x + (1-t) y } =F\brk{t x + (1-t) y}  \le tF(x)= tF^*(x)=tF^*(x) + (1-t)F^*(y).
\]

\end{proof}

 \begin{proof}[of \lemref{lem:short_convexity}]
 We prove the claim under the assumption that $F$ is convex  on $(1-\epsilon,1]$. The proof for the case of strict convexity is identical.
 Since $F|_{[1-\epsilon,1]}$ is (strictly) convex we have
 \beq
 \label{eq:extend_convex}
 F(x) \ge T_y(x):=F(y)+F_-'(y) (x-y) \, \, \, \text{ for every } \, \, x,y \in [1-\epsilon,1]. 
 \eeq

Let  $y \in [1-\epsilon,1]. $ Since $F$ is decreasing on $(0,1]$, $F_-'(y)  \le 0$, so $x \mapsto T_y(x)$ is decreasing. Since $T_y(1) \le F(1)= 0$, 
$T_y(x) \le 0 \le F(x)$ for every $x \ge 1$, so Inequality \eqref{eq:extend_convex} holds for every  $x \in [1-\epsilon,\infty)$ and every $y \in [1-\epsilon,1]. $ 

Let $x \in (0,1-\epsilon]$. For every $y <1$ sufficiently close to $1$, we have
\[
F(x) \ge F(1-\ep) \stackrel{(1)}{\ge} F(y)+|F_-'(y)|\ge F(y)+F_-'(y) (x-y)= T_y(x),
\]
where inequality $(1)$ follows from $\lim_{y \to 1^-} F_-'(y)=F_-'(1)=0$ together with $\lim_{y \to 1^-}F(y)=F(1)=0$. (We used the fact that the left derivative of a convex function is left-continuous.)
Thus, we proved that for $y <1$ sufficiently close to $1$, Inequality \eqref{eq:extend_convex} holds for every  $x \in (0,\infty)$, which means that $F=F|_{(0,\infty)}$ is convex at all such $y$.
 \end{proof}

\section{Proof of \lemref{lem:quantbound_vol_distortion_dist_SO22}}
\label{proofest1}

\lemref{lem:quantbound_vol_distortion_dist_SO22} is concerned with the distance of a matrix to the set $K$. We will therefore need the following claim: 
\begin{proposition}
\label{prop:dist_well_point_b}
Let $A \in M_2$ with $ \det A \ge 0$, and let $\sig_1\le \sig_2$ be its singular values. Then
\beq
\label{eq:quantbound_vol_distortion_dist_SO2_0}
\dist^2(A,K)=\begin{cases}
\frac{1}{2} \brk{\sig_1+\sig_2-1}^2,  & \text{ if  }\, \sig_2 \le \sig_1 + 1  \\
 \sigma_1^2+\brk{\sigma_2-1 }^2 , & \text{ if  }\,\sig_2 \ge \sig_1 + 1
\end{cases}
\eeq
\end{proposition}

We first use \propref{prop:dist_well_point_b} for proving \lemref{lem:quantbound_vol_distortion_dist_SO22}, then we prove it.
\begin{proof}[Of \lemref{lem:quantbound_vol_distortion_dist_SO22}]
Let $\sigma_1,\sigma_2$ be the singular values of $A$. Then 
\beq
\label{eq:quantbound_vol_distortion_dist_SO2_2}
\begin{split}
& \dist^2(A,\SOtwo)-(1-2\det A)= \\
& (\sigma_1-1)^2+(\sigma_2-1)^2-(1-2\sigma_1 \sigma_2) =  \\
& \sigma_1^2+\sigma_2^2+1-2\sigma_1-2\sigma_2+2\sigma_1 \sigma_2 = \\
&\brk{\sig_1+\sig_2-1}^2.
\end{split}
\eeq
The presence of the sum $\sig_1+\sig_2$ is no coincidence here! We represented the symmetric  polynomial $ P(\sigma_1,\sigma_2)=(\sigma_1-1)^2+(\sigma_2-1)^2$ as a polynomial in $\sigma_1 \sigma_2$ and $\sig_1+\sig_2$: $(x-1)^2+(y-1)^2=(1-2xy)+\brk{x+y-1}^2.$
\begin{comment}
Equation \eqref{eq:quantbound_vol_distortion_dist_SO2_2} implies that $\dist^2(A,\SOtwo) \ge 1-2\det A$ and equality holds exactly when $\sig_1+\sig_2=1$. This gives another proof for \lemref{lem:bound_vol_distortion_dist_SO22} in the regime where $\det A \le \bqrt$.
\end{comment}

Equation \eqref{eq:quantbound_vol_distortion_dist_SO2_2} and \propref{prop:dist_well_point_b} imply that  if $\sig_2 \le \sig_1 + 1$ then 
\[
\dist^2(A,\SOtwo)=(1-2\det A) +2\dist^2(A,K).
\]
Suppose that $\sig_2 \ge \sig_1 + 1$. Setting $x=\sig_1, y=\sig_2-1$, $x,y \ge 0$, hence 
\[
x^2+y^2 \le (x+y)^2 \le 2(x^2+y^2).
\] 
Equations \eqref{eq:quantbound_vol_distortion_dist_SO2_2} and
 \eqref{eq:quantbound_vol_distortion_dist_SO2_0} imply that $(x+y)^2=\dist^2(A,\SOtwo)-(1-2\det A)$ and $x^2+y^2=\dist^2(A,K)$ which completes the proof.
\end{proof}

\subsection{Computing $\dist(.,K)$ }
We prove \propref{prop:dist_well_point_b}. By  \defref{def:doublewell_combined} $K=\cup_{0\le s\le \frac{1}{4}}K_s$; we use the following lemma (that we prove below): 

\begin{lemma}
\label{lem:singlewell_dist_form}
Let $0\le \sigma_1 \le \sigma_2$, and let $A \in M_2$ satisfy $\det A \ge 0$. Then 
\beq
\label{eq:eq_doublewell_dist_form}
\dist^2(A,K_{\sig_1,\sig_2})=\brk{\sigma_1(A)-\sigma_1 }^2+\brk{\sigma_2(A)-\sigma_2 }^2.
\eeq
\end{lemma}


\begin{proof}[Of \propref{prop:dist_well_point_b}]

Define $F:[0,\frac{1}{4}] \to [0,\infty)$ by $F(s)=\dist^2(A,K_s)$. 
Since $K=\cup_{0\le s\le \frac{1}{4}}K_s$, 
\[
\dist^2(A,K)=\min_{0 \le s \le \frac{1}{4}}\dist^2(A,  K_s)=\min_{0 \le s \le \frac{1}{4}}F(s).
\]
We prove that  $\min_{0 < s \le \frac{1}{4}} F(s)$ equals the RHS of Equation \eqref{eq:quantbound_vol_distortion_dist_SO2_0}. 

Given $0 \le s \le \bqrt$, let $\sigma_1(s) \le \sigma_2(s)$ be the unique numbers satisfying 
\[
\sigma_1(s)\sigma_2(s)=s,\,\,\,\sigma_1(s)+\sigma_2(s)=1.
\] 
Since $K_s=K_{\sig_1(s),\sig_2(s)}$,  \lemref{lem:singlewell_dist_form} implies that 
\[
F(s)=\brk{\sigma_1-\sigma_1(s) }^2+\brk{\sigma_2-\sigma_2(s) }^2.
\]


Note that $\sigma_1(s)=\frac{1}{2} - \frac{\sqrt{1-4s}}{2},\sigma_2(s)=\frac{1}{2} + \frac{\sqrt{1-4s}}{2}$; thus $\sigma_i(s), F(s)$ are smooth functions of $s$ on $[0,\frac{1}{4})$ and continuous on $[0,\frac{1}{4}]$. We shall use the following lemma, which we prove at the end of the current proof. 
\begin{lemma}
\label{lem:doublewelldist} 
The function $F(s)=\dist^2(A,K_s)$ has a critical point $0\le s^* < \frac{1}{4}$ if and only if $\sig_1<\sig_2 \le \sig_1+1$. When such a critical point exists, it is unique and satisfies $F(s^*)=\frac{1}{2} \brk{\sig_1+\sig_2-1}^2$.
\end{lemma}

Next, we claim that if $\sig_1<\sig_2 \le \sig_1+1$, then $\min_{0 \le s \le \frac{1}{4}} F(s)=F(s^*)$. 
(One can prove that $F$ is convex, so any critical point is a global minimum, but we won't do that.)
The possible candidates for minimum points are interior critical points $s \in (0,\bqrt)$, and the endpoints $0,\bqrt$.  Thus, we need to show that 
\[
F(0) \ge \frac{1}{2} \brk{\sig_1+\sig_2-1}^2, F( \frac{1}{4}) \ge \frac{1}{2} \brk{\sig_1+\sig_2-1}^2.
\]
Since $\sig_1(\frac{1}{4})=\sig_2(\frac{1}{4})=\half$, $F( \frac{1}{4})=\sig_1^2+\sig_2^2+\half-\sig_1-\sig_2$; thus
\[
\begin{split}
& F( \frac{1}{4}) \ge \frac{1}{2} \brk{\sig_1+\sig_2-1}^2 \iff \\
& \sig_1^2+\sig_2^2+1-2\sig_1-2\sig_2+2\sig_1\sig_2 \le 2\sig_1^2+2\sig_2^2+1-2\sig_1-2\sig_2 \iff \\
& 0 \le \sig_1^2+\sig_2^2-2\sig_1\sig_2=(\sig_1-\sig_2)^2.
\end{split}
\]
Since $\sig_1(0)=0,\sig_2(0)=1$, $F(0)=\sig_1^2+\sig_2^2-2\sig_2+1.$ 
Thus, 
\[
\begin{split}
& F( 0) \ge \frac{1}{2} \brk{\sig_1+\sig_2-1}^2  \iff \\
&\sig_1^2+\sig_2^2+1-2\sig_1-2\sig_2+2\sig_1\sig_2 \le 2\sig_1^2+2\sig_2^2-4\sig_2+2 \iff \\
& -2\sig_1-2\sig_2+2\sig_1\sig_2 \le \sig_1^2+\sig_2^2-4\sig_2+1 \iff \\
& 2(\sig_2-\sig_1) \le (\sig_2-\sig_1)^2+1,
\end{split}
\]
which always holds since $2x \le x^2+1$ holds for every real $x$.

If $\sig_1=\sig_2$ or $\sig_2 > \sig_1 + 1$, then $F$ has no critical points. In these cases all is left to do is to compare $F(0)$ and $F( \bqrt)$. A direct computation shows that $F(0) \le F( \frac{1}{4})$ iff $\sig_2 \ge \sig_1 + 1/2$, from which the conclusion follows.
\end{proof}

\begin{proof}[Of \lemref{lem:doublewelldist}]

\[
2F'(s)=-\brk{\sigma_1-\sigma_1(s)}\sigma_1'(s)-\brk{\sigma_2-\sigma_2(s)}\sigma_2'(s).
\]
Since $\sigma_1(s)+\sigma_2(s)=1 \Rightarrow \sigma_1'(s)=-\sigma_2'(s)$, we get 
\[
2F'(s)=\sigma_1'(s)\brk{ \Delta \sigma- \Delta \sigma(s)},
\]
where $\Delta \sigma:=\sigma_2-\sigma_1, \Delta \sigma(s):=\sigma_2(s)-\sigma_1(s)$. 

$\sig_1'(s)>0$ so $F'(s)=0$ if and only if $\Delta \sigma= \Delta \sigma(s)$. Since $s \to \Delta \sig(s)=\sqrt{1-4s}$ is strictly decreasing, the uniqueness of the critical point is established. We now prove existence.
The condition $\sig_1 < \sig_2 \le \sig_1+1$  is necessary: $0 \le \sig_i(s) \le 1$ implies that  if $F'(s^*)=0$ then $\sigma_2-\sigma_1= \Delta \sigma(s^*) \le 1$. Furthermore, since $s^* < \frac{1}{4}$, $ \Delta \sigma(s^*)>0$, the equality $\Delta \sigma= \Delta \sigma(s^*)$ implies that $\sig_1<\sig_2$. 

To prove sufficiency, note that for every $0 < r \le 1$, there exists $s \in [0,\frac{1}{4})$ satisfying $\Delta \sigma(s)=r$. Now, let $s^*$ be the critical point. Since $\sig_1(s^*)+\sig_2(s^*)=1$, 
\[
-2\sig_1(s^*)=-1+\sig_2(s^*)-\sig_1(s^*)=\sig_2-\sig_1-1,
\]
thus
\[
\sig_1-\sig_1(s^*)=\frac{1}{2} \brk{2\sig_1-2\sig_1(s^*)}=\frac{1}{2} \brk{2\sig_1+\sig_2-\sig_1-1}=\frac{1}{2} \brk{\sig_1+\sig_2-1},
\]
so 
\[
F(s^*)=2\brk{\sigma_1-\sigma_1(s^*) }^2=\frac{1}{2} \brk{\sig_1+\sig_2-1}^2,
\]
where in the first equality we have used the implication $F'(s^*)=0 \Rightarrow \sigma_2-\sigma_2(s^*) = \sigma_1-\sigma_1(s^*)$.
This completes the proof.
\end{proof}


\begin{proof}\label{proof:lemsomething}[Of \lemref{lem:singlewell_dist_form}]

Given $X \in K_{\sig_1,\sig_2}$, we have 
\[
|A-X|^2=|A|^2 + |X|^2 -2\IP{A}{X}.
\] Since $|A|$ and $|X|=\sig_1^2+\sig_2^2$ are constant, we need to maximize $X \to \IP{A}{X}$ over $X \in K_{\sig_1,\sig_2}$.
By Von Neumann's trace inequality, 
\[
 \IP{A}{X}=\tr(A^TX) \le \sig_1(A)\sig_1+\sig_2(A)\sig_2.
 \]

It remains to show that this upper bound is realized by some $X \in K_{\sig_1,\sig_2}$. Using the bi-$\SOtwo$ invariance, we may assume that $A=\diag(\sig_1(A),\sig_2(A))$ is positive semidefinite and diagonal; taking $X=\diag(\sig_1,\sig_2)$ then realizes the bound. Thus
\[
\begin{split}
\dist^2(A,K_{\sig_1,\sig_2})&=\sig_1(A)^2+\sig_2(A)^2+\sig_1^2+\sig_2^2-2\sig_1(A)\sig_1-2\sig_2(A)\sig_2 \\
&=\brk{\sigma_1(A)-\sigma_1 }^2+\brk{\sigma_2(A)-\sigma_2 }^2.
\end{split}
\]
\begin{comment}
In the reduction of the problem to the diagonal positive semidefinite case, we explicitly use the assumption that $\det A \ge 0$. Indeed, let $A=U\Sigma V^T$ be the SVD of $A$. If $\det A>0$, then either both $U,V \in \SOtwo$ or both $U,V \in \Otwom$. In the latter case, we can multiply by $\diag\left( -1,1\right) $ from both sides of $\Sig$ to make them in $\SOtwo$. A similar argument works when $\det A=0$.
\end{comment}
\end{proof}
\section{Estimating $\dist(.,\CO)$}
\label{conformal_estimate}
\begin{proof}[Of \lemref{lem:quantbound_vol_distortion_dist_confSO22}]

First, expand
\beq
\label{eq:quantbound_vol_distortion_dist_SO2_2conf}
\begin{split}
& \dist^2(A,\SOtwo)-2(\sqrt{\det A}-1)^2= \\
& (\sigma_1-1)^2+(\sigma_2-1)^2-2(\sqrt{\sigma_1 \sigma_2}-1)^2 =  \\
& \brk{ \sigma_1^2+\sigma_2^2-2\sigma_1 \sigma_2 }-2\sigma_1-2\sigma_2+4\sqrt{\sigma_1 \sigma_2} = \\
&\brk{\sig_2-\sig_1}^2-2\brk{\sqrt \sig_2-\sqrt \sig_1}^2 \le \brk{\sig_2-\sig_1}^2=2\dist^2(A,\CO).
\end{split}
\eeq
Now, 
\beq
\label{eq:quantbound_vol_distortion_dist_SO2_2confaa}
\brk{\sig_2-\sig_1}^2-2\brk{\sqrt \sig_2-\sqrt \sig_1}^2 = 
\brk{\sqrt \sig_2-\sqrt \sig_1}^2 \brk{\brk{\sqrt \sig_2+\sqrt \sig_1}^2-2}.
\eeq
By the AM-GM inequality, if $\det A \ge \bqrt$ then 
\[
\brk{\sqrt \sig_2+\sqrt \sig_1}^2 \ge 4\sqrt{  \sig_1\sig_2} \ge 2,
\] 
which implies 
\beq
\label{eq:quantbound_vol_distortion_dist_SO2_2confbb}
\brk{\sqrt \sig_2+\sqrt \sig_1}^2-2 \ge \brk{\sqrt \sig_2+\sqrt \sig_1}^2 - 4\sqrt{  \sig_1\sig_2}= \brk{\sqrt \sig_2-\sqrt \sig_1}^2.
\eeq
Equations \eqref{eq:quantbound_vol_distortion_dist_SO2_2confaa}, \eqref{eq:quantbound_vol_distortion_dist_SO2_2confbb} complete the proof.
\end{proof}

\section{The Euler-Lagrange equation of $\dist^2(.,\SO)$}
\label{sec:EL_deriv}
In this section we prove that the Euler-Lagrange equation of the functional of $E_2$ is
\[
\delta \brk{d\phi-O(d\phi)}=0,
\]
where $O(d\phi)$ is the orthogonal polar factor of $d\phi$. (The derivation of the EL equations for $p \neq 2$ follows from the special case of $p=2$.)

For brevity, we show the derivation only for the Euclidean case where $\M=\Omega \subseteq \R^n, \N=\R^n$ are endowed with the usual flat metrics; the general Riemannian case follows in a similar fashion.

Let $\GLp$ be the group of real $n \times n$ matrices having positive determinant, and let $O:\GLp \to \SOn$ map $A\in \GLp$ into its \emph{orthogonal polar factor}, i.e.
\[
O(A)=A\brk{\sqrt{A^TA}}^{-1}.
\]
$\sqrt{A^TA}$ denotes the unique symmetric positive-definite square root of $A^TA$. The map $O$ is smooth. We use the following observation:

\begin{lemma}
\label{lem:orthogonal_polar_deriv_orthogonal_relat}
Let $A \in \GLp$. Given $B\in M_n$ write $\dot O=dO_A(B)$. Then for every $B \in M_n$,
\[
\IP{\dot O}{O}=\IP{\dot O}{A}=0.
\]
\end{lemma}

\begin{proof}
The equality $\IP{\dot O}{O}=0$ follows from differentiating $\IP{O}{O}=n$. Now,
\[
\dot O \in T_{O}\SOn=OT_{\id}\SOn=O\skew
\]
implies that $ \dot O=OS$ for some $S \in \skew$. Thus,
\[
\IP{\dot O}{A} = \IP{OS}{OP} = \IP{S}{P} =0,
\]
where the last equality follows from the fact that the spaces of symmetric matrices and skew-symmetric matrices are orthogonal.
\end{proof}

\paragraph{Derivation of the EL equation}

Recall that
\[
E_2(\phi)=\int_{\Omega} \dist^2\brk{d\phi,\SOn}=\int_{\Omega}|d\phi-O(d\phi)|^2.
\] 
Let $\phi \in C^2(\Omega,\R^n)$, and let $\phi_t=\phi+tV$ for some $C^2$ vector field $V:\Omega \to \R^n$. Then, 
\beq 
\begin{split}
\label{eq:energy_time_deriv_Riemann}
\half \left. \deriv{}{t} E\left( \phi_t \right) \right|_{t=0}&= \int_{\Omega}  \left. \IP{ \covder_{\pd{}{t}} d\phi_t  - \covder_{\pd{}{t}} O(d\phi_t)   }{ d\phi_t-O(d\phi_t)}  \right|_{t=0} 
 dx \\
 & =\int_{\Omega} \IP{ \left. \covder_{ \pd{}{t} } d\phi_t     \right|_{t=0}  }{ d\phi-Q(d\phi)} dx=\int_{\Omega} \IP{ \nabla V  }{ d\phi-Q(d\phi)} dx \\
 &=-\int_{\Omega}  \IP{  V  }{ \div\brk{d\phi-Q(d\phi)}}  dx.
 \end{split}
\eeq
The passage from the first to the second line relied upon the fact that $\covder_{\pd{}{t}} O(d\phi_t)$ is orthogonal to $d\phi_t,O(d\phi_t) $, which essentially follows from \lemref{lem:orthogonal_polar_deriv_orthogonal_relat}.


\section{Additional proofs}
\label{sec:add_proof}
\begin{lemma}
\label{lem:minfunc_prop}
The function $F$ defined in \eqref{eq:min_uniform0} is well-defined and continuous; 
it is strictly decreasing on $(0,1]$ and strictly increasing on $[1,\infty)$.  
\end{lemma}

\begin{proof}
First, suppose that $s \le 1$; then the minimum is obtained at a point $(a,b)$ where both $a,b \le 1$. Indeed, if $a>1$ (and so $b <s \le 1$), we can replace $a$ by $1$ and $b$ by $s$ to get the same product with both numbers closer to $1$. Thus, it suffices to show that the minimum exists when $0<a,b \le 1$. Now, $b \le 1 \Rightarrow s=ab \le a$, and similarly $b \ge s$. So, the problem reduces to proving existence of a minimum over the compact set $\{ (a,b) \in [s,1]^2 \, | \, ab=s\}$. Since $f$ was assumed continuous we are done.

Next, we prove that $F$ is strictly decreasing on $(0,1]$. Indeed, let $0< s_1 < s_2 \le 1$, and suppose that $F(s_1)=f(a)+f(b)$, for some $(a,b) \in [s_1,1]^2, ab=s_1$. Choose a smooth path $(a(t),b(t))$ from $(a,b)$ to $(1,1)$, where $a(t),b(t)$ are both strictly increasing. Then for $t>0$
\[
F\brk{a(t)b(t)} \le f(a(t))+f(b(t)) < f(a)+f(b)=F(s_1).
\]
Since $a(t)b(t)$ movies continuously from $s_1$ to $1$, it hits $s_2$ at some time $t>0$, which establishes the claim. Finally, a symmetric argument shows that if $s \ge 1$, then the minimum is obtained in  $\{ (a,b) \in [1,s]^2 \, | \, ab=s\}$, and that $F$ is strictly increasing on $[1,\infty)$. Proving $F$ is continuous is routine and we omit it.
\end{proof}

\begin{proof}[of \propref{prop:div_inf_min}]
Suppose that $(\sqrt s_n,\sqrt s_n)$ is a minimizer of \eqref{eq:min_uniform0} for some sequence $s_n \in (0,1)$ which converges to zero.  We prove that $f$ diverges to $\infty$ at zero. Recasting everything in terms of $g(x)=f(e^x)$, we get  \[
g(\frac{x + y}{2}) \le \frac{g(x) + g(y)}{2}  \, \,\, \, \text{ whenever } \, \, x,y \le 0 \, \, \, \text{ and  } x+y=\lambda_n,
\]
where $\lambda_n=\log s_n \to -\infty$. Choose a subsequence $(\lambda_{n_k})$ such that $
 \lambda_{n_k} < 2 \lambda_{n_{k-1}}$ for all $k$. Choosing $x=0$ and $y= \lambda_{n_k}$ in the condition above and the monotonicity of $g$ give
\[
 g(\lambda_{n_{k-1}}) \le  g\left( \frac 12 \lambda_{n_k}\right) \le \frac 12 g(\lambda_{n_k})
\]
so
\[
 g(\lambda_{n_k}) \ge 2 g(\lambda_{n_{k-1}}) \ge 2^2 g(\lambda_{n_{k-2}})
\ge \dots \ge 2^{k-1}  g(\lambda_{n_{1}}) \, .
\]
 \end{proof}

\begin{proof}[Of \lemref{lem:bound_vol_distortion_dist_SO22}]
Since the problem is bi-$\SOtwo$-invariant, using SVD we can assume that $A=\diag(a,b)$ is diagonal with $a,b>0$. We need to compute 
\[
 F(s)=\min_{a,b \in \mathbb{R}^+,ab=s} (a-1)^2+(b-1)^2. 
 \]
Using Lagrange's multiplier, there exists $\lambda$ such that $\brk{a-1,b-1}=\lambda(b,a)$.
Thus $a(1-a)=b(1-b)$ which implies $a=b$ or $a=1-b$.
%
%
%
%
In the latter case $s=ab=b(1-b)$. Since $a=1-b,b,s$ are positive, we must have $0<b<1,0<s\le\frac{1}{4}$. (since $\max_{0<b<1} b(1-b)=\frac{1}{4}$).
We then have 
\[
(a-1)^2+(b-1)^2 =b^2+(b-1)^2=2b(b-1)+1=1-2s.
\]
\begin{itemize}
\item If $s \ge \bqrt$ then there is only one critical point $(a,b)=(\sqrt s,\sqrt s)$, so the minimum is obtained exactly when $a=b$ and $F(s)=2(\sqrt{s}-1)^2.$ 

\item If $s \le \bqrt$, there are up to $3$ to critical points: $(\sqrt s, \sqrt s),(x,1-x),(1-x,x)$ when $0<x<1$ satisfies $x(1-x)=s$. (For $s=\bqrt $ they all merge into a single point. For $s<\bqrt$ these are $3$ distinct points.) 

To decide which point is the global minimizer, we need to compare the values of the objective function at the critical points, which are $2(\sqrt{s}-1)^2,1-2s.$
Since 
\[
1-2s \le 2(\sqrt{s}-1)^2 \iff (2\sqrt{s}-1)^2 \ge 0,
\] 
$F(s)=1-2s$ for $s \le \frac{1}{4}$. 
\end{itemize}
\end{proof}

\begin{proof}[of \lemref{lem:Piola_approach_limit}]
By \propref{prop:geometric_cof_char_well}, 
\[
d\phi_x-O(d\phi_x)=\al \Cof d\phi_x
\]
if and only if $d\phi_x \in K$ and $\al=-1$ or $d\phi_x$ is conformal and $\sig_i(d\phi_x)=\frac{1}{1-\al}$. In both cases $\sig_1(d\phi_x)+\sig_2(d\phi_x)=\frac{2}{1-\al}$ is constant. If $\al=-1$, then $d\phi \in K$, and if $\al \neq -1$ then $d\phi \in \CO$ with $\sig_i(d\phi)=\frac{1}{1-\al}$, i.e.~$\phi$ is a homothety.

Now, define $H(\phi)=\dist^{p-2}\brk{d\phi,\SO(\g,\phi^*\h)} $. Let $p \neq 2$ and suppose that 
\[
H(\phi)\brk{d\phi-O(d\phi)}=\al \Cof d\phi.
\]
Since we assumed $J\phi>0$, $\Cof d\phi$ is invertible, hence $H(\phi) \neq 0$, and  
\[
d\phi-O(d\phi)=\frac{\al}{H(\phi)} \Cof d\phi.
\]

 \propref{prop:geometric_cof_char_well} implies that either $H(\phi)(x)=-\al$ and $d\phi_x \in K$, or that $d\phi_x$ is conformal and 
 \[
 \sig(d\phi_x)=\frac{1}{1-\al/H(\phi)(x)}=\frac{1}{1-\frac{\al}{(\sqrt{2}|\sig(d\phi_x)-1|)^{p-2}}}.
 \]
 Thus $ \sig(d\phi_x)$ is a solution for the equation
  \[
 \sig=\frac{1}{1-\frac{\be}{|\sig-1|^{p-2}}},
 \]
 which has a finite number of solutions $\{\al_1,\dots,\al_k\}$.
 
 We showed that for every $x\in \M$, $d\phi_x \in K$, or $d\phi_x$ is conformal with $ \sig(d\phi_x) \in \{\al_1,\dots,\al_k\}$.
 Thus $x \mapsto \sig_1(d\phi_x)+\sig_2(d\phi_x)$ is a continuous function on $\M$, which takes  values in a finite set $\{1,2\al_1,\dots,2\al_k\}$, hence it must be constant.
Thus, $d\phi \in K$ or $\phi$ is a homothety with $\sig(d\phi) \in \{\al_1,\dots,\al_k\}$. If $d\phi \in K$, then $H(\phi)=-\al$ is constant. Thus 
\[
 \sig_1(d\phi)+ \sig_2(d\phi)=1,  \brk{\sig_1(d\phi)-1}^2+\brk{\sig_2(d\phi)-1}^2
\]
are constants, which implies that $ \sig_1(d\phi), \sig_2(d\phi)$ are constants as required.
 
\end{proof}


\bibliographystyle{unsrt}

\begin{thebibliography}{10}

\bibitem{ESK08}
E.~Efrati, E.~Sharon, and R.~Kupferman.
\newblock Elastic theory of unconstrained non-{Euclidean} plates.
\newblock {\em Journal of the Mechanics and Physics of Solids}, 57:762--775,
  2009.

\bibitem{KES07}
Y.~Klein, E.~Efrati, and E.~Sharon.
\newblock Shaping of elastic sheets by prescription of non-{Euclidean} metrics.
\newblock {\em Science}, 315:1116 -- 1120, 2007.

\bibitem{KVS11}
Y.~Klein, S.~Venkataramani, and E.~Sharon.
\newblock Experimental study of shape transitions and energy scaling in thin
  non-euclidean plates.
\newblock {\em PRL}, 106:118303, 2011.

\bibitem{DCGLL13}
A.~Danescu, C.~Chevalier, G.~Grenet, Ph. Regreny, X.~Letartre, and J.L.
  Leclercq.
\newblock Spherical curves design for micro-origami using intrinsic stress
  relaxation.
\newblock {\em Applied Physics Letters}, 102(12):123111, 2013.

\bibitem{AKMMS16}
H.~Aharoni, J.~Kolinski, M.~Moshe, I.~Meirzada, and E.~Sharon.
\newblock Internal stresses lead to net forces and torques on extended elastic
  bodies.
\newblock {\em Physical review letters}, 117(12):124101, 2016.

\bibitem{Sil01}
M.~{\v{S}}ilhav\'{y}.
\newblock Rank-1 convex hulls of isotropic functions in dimension 2 by 2.
\newblock {\em Proceedings of Partial Differential Equations and Applications
  (Olomouc, 1999)}, 126:521--529, 2001.

\bibitem{Dol12}
G.~Dolzmann.
\newblock Regularity of minimizers in nonlinear elasticity -- the case of a
  one-well problem in nonlinear elasticity.
\newblock {\em TECHNISCHE MECHANIK}, 32:189--194, 2012.

\bibitem{Har58}
P.~Hartman.
\newblock On isometries and on a theorem of liouville.
\newblock {\em Mathematische Zeitschrift}, 69:202--210, 1958.

\bibitem{CH70}
E.~Calabi and P.Hartman.
\newblock On the smoothness of isometries.
\newblock {\em Duke Math. J.}, 37(4):741--750, 12 1970.

\bibitem{Tay06}
M.~Taylor.
\newblock Existence and regularity of isometries.
\newblock {\em Transactions of the American Mathematical Society},
  358(6):2415--2423, 2006.

\bibitem{kupferman2019reshetnyak}
R.~Kupferman, C.~Maor, and A.~Shachar.
\newblock Reshetnyak rigidity for {Riemannian} manifolds.
\newblock {\em Archive for Rational Mechanics and Analysis}, 231(1):367--408,
  2019.

\bibitem{351550}
R.~Bryant (https://mathoverflow.net/users/13972/robert bryant).
\newblock Are all maps $\mathbb{R}^2 \to \mathbb{R}^2$ with fixed singular
  values affine?
\newblock MathOverflow.
\newblock URL:https://mathoverflow.net/q/351550 (version: 2020-02-01).

\bibitem{RobObs}
R.~Bryant.
\newblock Communication on the mathoverflow website (2020), available online at
  https://mathoverflow.net/questions/376018/metric-obstructions-for-area-preserving-diffeomorphisms-with-constant-singular-v.

\bibitem{375931}
R.~Bryant (https://mathoverflow.net/users/13972/robert bryant).
\newblock A diffeomorphism of the torus with constant singular values.
\newblock MathOverflow.
\newblock URL:https://mathoverflow.net/q/375931 (version: 2020-11-08).

\bibitem{bertoldi2017flexible}
K.~Bertoldi, V.~Vitelli, J.~Christensen, and M.~van Hecke.
\newblock Flexible mechanical metamaterials.
\newblock {\em Nature Reviews Materials}, 2(11):1--11, 2017.

\bibitem{stoop2015curvature}
N.~Stoop, R.~Lagrange, D.~Terwagne, P.M. Reis, and J.~Dunkel.
\newblock Curvature-induced symmetry breaking determines elastic surface
  patterns.
\newblock {\em Nature materials}, 14(3):337--342, 2015.

\bibitem{Res67}
Yu.~G. Reshetnyak.
\newblock On the stability of conformal mappings in multidimensional spaces.
\newblock {\em Sibirskii Matematicheskii Zhurnal}, 8(1):91--114,
  January--February 1967.

\bibitem{Mul90}
S.~M{\"{u}}ller.
\newblock Higher integrability of determinants and weak convergence in {$L^1$}.
\newblock {\em Journal f\''ur die reine und angewandte Mathematik}, 412:20--34,
  1990.

\bibitem{giaquinta1998cartesian1}
M.~Giaquinta, G.~Modica, and J.~Soucek.
\newblock {\em Cartesian Currents in the Calculus of Variations II: Variational
  Integrals}, volume~1.
\newblock Springer Science \& Business Media, 1998.

\bibitem{381194}
A.~Shachar (https://mathoverflow.net/users/46290/asaf shachar).
\newblock Does weak continuity of {Jacobians} hold for non nondegenerate maps?
\newblock MathOverflow.
\newblock URL:https://mathoverflow.net/q/381194 (version: 2021-01-15).

\bibitem{Cia88}
P.~G. Ciarlet.
\newblock {\em Mathematical Elasticity, Volume 1: Three-dimensional
  elasticity}.
\newblock Elsevier, 1988.

\bibitem{sivaloganathan2010global}
J.~Sivaloganathan and S.~J. Spector.
\newblock On the global stability of two-dimensional, incompressible, elastic
  bars in uniaxial extension.
\newblock {\em Proceedings of the Royal Society A: Mathematical, Physical and
  Engineering Sciences}, 466(2116):1167--1176, 2010.

\bibitem{mora2010explicit}
C.~Mora-Corral.
\newblock Explicit energy-minimizers of incompressible elastic brittle bars
  under uniaxial extension.
\newblock {\em Comptes Rendus Mathematique}, 348(17-18):1045--1048, 2010.

\bibitem{hajlasz1999sobolev}
P.~Haj{\l}asz.
\newblock Sobolev mappings, co-area formula and related topics.
\newblock 1999.

\bibitem{bony2010square}
J-M. Bony, F.~Colombini, and L.~Pernazza.
\newblock On square roots of class {$C^m$} of nonnegative functions of one
  variable.
\newblock {\em Annali della Scuola Normale Superiore di Pisa-Classe di
  Scienze}, 9(3):635--644, 2010.

\bibitem{374383}
P.~Hajlasz (https://mathoverflow.net/users/121665/piotr hajlasz).
\newblock Is {$L^1$} strong convergence of {Jacobians} valid for maps between
  manifolds?
\newblock MathOverflow.
\newblock URL:https://mathoverflow.net/q/374383 (version: 2020-10-20).

\bibitem{Lee13}
J.~M. Lee.
\newblock {\em Introduction to Smooth Manifolds}.
\newblock Springer, 2nd edition, 2013.

\bibitem{253994}
Mizar (https://mathoverflow.net/users/36952/mizar).
\newblock Are metric isometries smooth at the boundary?
\newblock MathOverflow.
\newblock URL:https://mathoverflow.net/q/253994 (version: 2016-11-05).

\bibitem{EL83}
J.~Eells and L.~Lemaire.
\newblock {\em Selected topics in harmonic maps}, volume~50.
\newblock American Mathematical Soc., 1983.

\bibitem{KS19}
R.~Kupferman and A.~Shachar.
\newblock A geometric perspective on the piola identity in riemannian settings.
\newblock {\em Journal of Geometric Mechanics}, 11(1):59--76, 2019.

\bibitem{Eva98}
L.~C. Evans.
\newblock {\em Partial Differential Equations}.
\newblock American Mathematical Society, 1998.

\bibitem{383251}
R.~Bryant (https://mathoverflow.net/users/13972/robert bryant).
\newblock Local obstructions for maps with constant singular values.
\newblock MathOverflow.
\newblock URL:https://mathoverflow.net/q/383251 (version: 2021-02-08).

\bibitem{desimone2002macroscopic}
A.~DeSimone and G.~Dolzmann.
\newblock Macroscopic response of nematic elastomers via relaxation of a class
  of so (3)-invariant energies.
\newblock {\em Archive for rational mechanics and analysis}, 161(3):181--204,
  2002.

\bibitem{3163368}
Dap (https://math.stackexchange.com/users/467147/dap).
\newblock Can we choose smoothly the singular vectors of a matrix?
\newblock Mathematics Stack Exchange.
\newblock URL:https://math.stackexchange.com/q/3163368 (version: 2019-03-31).

\bibitem{opt20}
I.~Pinelis (https://mathoverflow.net/users/36721/iosif pinelis).
\newblock Is the optimum of this problem convex in the constraint parameter.
\newblock MathOverflow.
\newblock URL:https://mathoverflow.net/q/357467 (version: 2020-04-14).

\bibitem{neff2014logarithmic}
P.~Neff, Y.~Nakatsukasa, and A.~Fischle.
\newblock A logarithmic minimization property of the unitary polar factor in
  the spectral and frobenius norms.
\newblock {\em SIAM Journal on Matrix Analysis and Applications},
  35(3):1132--1154, 2014.

\bibitem{neff2016geometry}
P.~Neff, B.~Eidel, and R.J. Martin.
\newblock Geometry of logarithmic strain measures in solid mechanics.
\newblock {\em Archive for Rational Mechanics and Analysis}, 222(2):507--572,
  2016.

\bibitem{lankeit2014minimization}
J.~Lankeit, P.~Neff, and Y.~Nakatsukasa.
\newblock The minimization of matrix logarithms: On a fundamental property of
  the unitary polar factor.
\newblock {\em Linear Algebra and its Applications}, 449:28--42, 2014.

\bibitem{KS16}
R.~Kupferman and A.~Shachar.
\newblock On strain measures and the geodesic distance to {$SO_n$} in the
  general linear group.
\newblock {\em Journal of Geometric Mechanics}, 8(4):437--460, 2016.

\end{thebibliography}
\end{document}